\documentclass[a4paper,12pt]{article}
\usepackage{amssymb}
\usepackage{latexsym}
\usepackage{amsmath}
\usepackage{amsthm}
\usepackage{amsfonts}
\usepackage{amssymb}
\usepackage{amsmath,amscd}
\usepackage{graphicx}
\usepackage{color}

\definecolor{babyblue}{rgb}{0.54, 0.81, 0.94}
\definecolor{bittersweet}{rgb}{1.0, 0.44, 0.37}
\definecolor{brightmaroon}{rgb}{0.76, 0.13, 0.28}
\definecolor{icterine}{rgb}{0.99, 0.97, 0.37}
\definecolor{indiagreen}{rgb}{0.07, 0.53, 0.03}
\definecolor{aquamarine}{rgb}{0.5, 1.0, 0.83}
\definecolor{coolblack}{rgb}{0.0, 0.18, 0.39}
   	\definecolor{cobalt}{rgb}{0.0, 0.28, 0.67}
	\definecolor{amber}{rgb}{1.0, 0.75, 0.0}
	\definecolor{azure(colorwheel)}{rgb}{0.0, 0.5, 1.0}
	 \definecolor{shamrockgreen}{rgb}{0.0, 0.62, 0.38}
	 \definecolor{almond}{rgb}{0.94, 0.87, 0.8}
\definecolor{arylideyellow}{rgb}{0.91, 0.84, 0.42}
\definecolor{bubblegum}{rgb}{0.99, 0.76, 0.8}
\definecolor{babypink}{rgb}{0.96, 0.76, 0.76}
 	\definecolor{ballblue}{rgb}{0.13, 0.67, 0.8}

\newtheorem{Th}{Theorem}
\newtheorem{Prop}{Proposition}
\newtheorem{Co}{Corollary}
\newtheorem{Lm}{Lemma}

\newtheorem{Rm}{Remark}

\newcommand{\be}{\begin{equation}}
\newcommand{\ee}{\end{equation}}
\newcommand{\bes}{\begin{equation*}}
\newcommand{\ees}{\end{equation*}}
\newcommand{\brac}[1]{\left (#1 \right )}
\def\11{1\!\!1} 
\newcommand{\K}{\mathbb{H}}
\newcommand{\R}{\mathbb{R}}

\newcommand{\C}{\mathbb{C}}
\newcommand{\Z}{\mathbb{Z}}

\def\XXint#1#2#3{{\setbox0=\hbox{$#1{#2#3}{\int}$ }
\vcenter{\hbox{$#2#3$ }}\kern-.6\wd0}}

\def\ti{\tilde}
\def\lf{\left}
\def\rg{\right}

\def\al{\alpha}
\def\la{\lambda}

\def\ep{\varepsilon}
\def\eps{\varepsilon}
\def\ds{\displaystyle}
\def\ov{\overline}
\def\Om{\Omega}
\def\om{\omega}
\def\p{\partial}

\def\nablap{\nabla^\perp}

\AtBeginDocument{
  \let\div\relax
  \DeclareMathOperator{\div}{div}
}

\begin{document}

\title{ Critical Chirality in Elliptic Systems}

\author{ Francesca Da Lio and Tristan Rivi\`ere\footnote{Department of Mathematics, ETH Zentrum,
CH-8093 Z\"urich, Switzerland.}}

%\date{ }
\maketitle
{\bf Abstract :}{\it We establish the regularity in $2$ dimension of $L^2$ solutions to critical elliptic systems in divergence form involving chirality operators of finite $W^{1,2}$-energy.}

\par

\section{Introduction}
In \cite{Riv1} the second author discovered a {\it compensation phenomenon} for the linear elliptic systems of the form
\be
\label{intro-1}
-\Delta u=\Omega\cdot \nabla u\quad\mbox{ in } D^2\,,
\ee
where $u\in W^{1,2}(D^2,{\R}^n)$, $D^2=B^2(0,1)$ denotes the open unit ball centered at $(0,0)$ and $\Omega$ is an $L^2$ map into  the antisymmetric matrices of ${\R}^2$ vectors. That is to say there exists a matrix $(\Om_i^j)_{i,j=1\cdots n}$ of $L^2$ functions into ${\R}^2$ such that
\[
\forall\ i=1\cdots n\quad\quad-\Delta u_i=\sum_{j=1}^n\Om_i^j\cdot\nabla u_j\quad\quad\mbox{ and }\quad \Om_i^j=-\Om_j^i\quad\forall \ i,j=1\cdots n\,.
\]A-priori the system
(\ref{intro-1}) is critical for the chosen norms, with a right hand side in $L^1$. Without the anti-symmetry of $\Omega$ no improved regularity has to be expected in general, while $W^{1,2}$ solutions to (\ref{intro-1})  for $\Om\in L^2({D^2}, {\R}^2\otimes so(n))$ are known to be in $\bigcap_{p<2}W^{2,p}_{loc}(D^2)$.

\medskip

One of the main strategy introduced in \cite{Riv1} was to use the antisymmetry of $\Omega$ in order to construct a ``gauge'' $A\in L^\infty\cap W^{1,2}(D^2,Gl_n({\R}))$ satisfying
\[
\mbox{div}(\nabla_\Omega A):=\mbox{div}(\nabla A- A\,\Omega)=0
\]
Taking a ``primitive'' $\nabla^\perp B=(-\p_{x_2}B, \p_{x_1}B):=\nabla_\Omega A\in L^2(D^2,M_n({\R})\otimes {\R}^2)$ the system (\ref{intro-1}) becomes equivalent to the conservation law
\be
\label{intro-2}
\mbox{div}(A\nabla u)=\nabla^\perp B\cdot\nabla u
\ee
The Jacobian form of the right-hand-side of (\ref{intro-2}) permits to use now classical integrability by compensation phenomena originally discovered
by H.Wente \cite{Wen} and related to the ones by R.Coifman, Rochberg and Weiss \cite{CRW} (see also \cite{CLMS}).

Following the main ideas  of \cite{Riv1}, extensions of this compensation phenomenon were obtained in \cite{Riv2} for critical systems of the  form (for $m>2$)
\be
\label{intro-3}
\Delta v=\Omega\, v\quad\mbox{ in } B^m\,,
\ee
where $B^m$ denote the $m$-dimensional ball centered at $0$ and with radius $1$, $\Om\in L^{m/2}(B^m,so(n))$ and $v\in L^{m/(m-2)}(B^m,{\R}^n)$ as well as for systems of the form 
\be
\label{intro-4}
(-\Delta)^{1/4} v=\Om\, v\quad\mbox{ in }{\R}
\ee
where this time $v\in L^2({\R},{\R}^n)$ and $\Omega\in L^2({\R},so(n))$ (see \cite{DLR2}). More recently the two authors are extending their results to non local right-hand-side of the form
\be
\label{intro-5}
(-\Delta)^{1/4}v=\int_{\R}H(x,y)\, v(y)\ dy
\ee
where pointwise antisymmetry has to be replaced by the more general notion of anti-self-duality of the underlying non-local operator $K(x,y)$ where
$K(x,y):=H(x,y)-\om(x)\,\delta_{x=y}\in L^1_{loc}({\R}^2)$ (see \cite{DLR3}).

In the present work we are exhibiting a new compensation phenomenon which does not enter in none of the previous existing ones. Our main result is the following
\begin{Th}
\label{th-intro-1}
Let $S\in \dot{W}^{1,2}(\R^2,O(n))$, \footnote{ $O(n)$ denotes the  group of orthogonal $n\times n$ matrices, $SO(n)$ is the group of orthogonal $n\times n$ matrices with determinant $1$. $U(n)$ is the group of unitary $n\times n$ matrices and $SU(n)$ is the Lie group of $n\times n$ unitary matrices with determinant $1$. The Lie algebra of $U(n)$ consists of $n\times n$ skew-Hermitian matrices, with the Lie bracket given by the commutator.}  such that $S^2=id_n$ and let
$u\in L^2(\R^2,{\R}^n)$ be a solution of the following linear elliptic system in divergence form
\begin{eqnarray}
\label{intro-6}
\mbox{div}\,(S\,\nabla u)&=&\sum_{j=1}^n\mbox{div}\,(S_{ij}\,\nabla u^j)=\sum_{j=1}^n\sum_{\alpha=1}^2\frac{\partial}{\partial_{x_{\alpha}}}\left(S_{ij}u_{x_{\alpha}}^j\right)=0
\end{eqnarray}
Then $\ds u\in \bigcap_{p<2}W^{1,p}_{loc}(\R^2,{\R}^n)$.\hfill $\Box$
\end{Th}
\begin{Rm}
\label{rm-intro-0}
The system (\ref{intro-6}) is elliptic with principal symbol $|\xi|^2 S$. It is however not {\it strongly elliptic} in the sense of Legendre Hadamard \footnote{We recall that
a matrix of coefficients $(A^{\alpha,\beta}_{ij})^{1\le\alpha,\beta\le m}_{1\le i,j\le n}$ satisfies the {\em strong ellipticity condition}, or the {\em Legendre-Hadamard condition} if there exists $\lambda>0$ such that 
$$A^{\alpha,\beta}_{ij}\xi_\alpha\xi_\beta\eta^i\eta^j\ge\lambda |\xi|^2|\eta|^2,~~\mbox{for all $\xi\in\R^m,~\eta\in\R^n.$}$$
In the case of the system \eqref{intro-6} the matrix of coefficients is given by $A^{\alpha,\beta}_{ij}:=S_{ij}\delta_{\alpha,\beta}$ where
$1\le \alpha,\beta\le 2$, $1\le i,j\le n$ and $\delta_{\alpha,\beta}$ denotes the Kronecker's operator.}
since obviously $<S\la,\la>$ can change sign as $\la$ varies.\hfill $\Box$
\end{Rm}
{\begin{Rm}
\label{rm-intro-1}
  Structural conditions on $S$ for the regularity are necessary in the following sense. In \cite{JMS} an $L^2$ solution to
\[
\mbox{div}\,(A\,\nabla u)=0
\]  
is produced where $A\in W^{1,2}(D^2,\mbox{Sym}(2))$ and 
$A$ is satisfying the strong ellipticity condition\footnote{The matrix $A$ is acting on the different vertical components of $\nabla u$, $u$ is in fact scalar in this case while in theorem~\ref{th-intro-1} the matrix $S$ acts on the horizontal components of $\nabla u$ that is $\nabla u_1,\cdots,\nabla u_n$.} $$<A(x)\xi,\xi>\simeq |\xi|^2$$ uniformly on $D^2$ but $u\notin W^{1,p}_{loc}(D^2,{\R})$ for any $p>1$. \par
We also observe that we cannot expect in Theorem \ref{th-intro-1} that $u\in W^{1,2}_{loc}(D^2)$. Actually if we set $w=Su$, $w$ solves
$\Delta w=\div(\nabla S Sw)$. Such a PDE bootstraps in $W^{1,p}_{loc}(D^2)$ for $p<2$ but not in $W^{1,2}_{loc}(D^2)$.
If $w\in W^{1,p}_{loc}(D^2)$  for $p<2$ then $w\in L^{p^*}(D^2)$ with $p^*=\frac{2p}{2-p}$. By injecting such an information into the equation we get
that $\nabla S Sw\in L^p_{loc}(D^2)$ (since $p=\frac{2p^*}{p^*+2}$) and therefore we come back to the initial information that $\nabla w \in L^p_{loc}(D^2)$. This is not the case if $w\in W^{1,2}_{loc}(D^2)$. This  would actually  imply that $w\in L^{q}_{loc}(D^2)$ for every $q\ge 2$ and from the equation we deduce  that
$\nabla w\in L^{\frac{2q}{q+2}}_{loc}$ which is a lost of information from the initial one since $1\le \frac{2q}{q+2}<2$.
 
\end{Rm}}
\begin{Rm}
\label{rm-intro-2} Contrary  to the case of the  systems (\ref{intro-1}) in \cite{Riv1}, we have not found yet striking applications in geometry or physics of systems (\ref{intro-6})
while nevertheless they look very ``natural'' and enjoy numerous formulations  that we are going to present in this work. The system (\ref{intro-6}) is nothing but the Harmonic Map Equation into a pseudo-riemannian manifold (see remark~\ref{rm-pseudo}). The formulation  using Dirac operator below (see \ref{intro-10}) moreover corresponds to the {\bf Weierstrass representation} of {\bf Lagrangian surfaces} in
four-dimensional space by H\'elein and Romon (\cite{HR} Theorem 1). The assumption $u\in L^2$ is  also faithful to the Hilbert Space framework in mathematical physics\footnote{Original works in mathematical physics which have nourish the growth of analysis with problems from quantum mechanics, such as the study of Schr\"odinger semigroups for instance \cite{Sim}...etc, take the $L^2$
space and not the ``energy space'' $W^{1,2}$ as the ''configuration space''.} and makes this function space natural in that sense.\hfill $\Box$
\end{Rm}
Behind the proof of theorem~\ref{th-intro-1} there is an $\epsilon-$regularity type of estimate which implies the following concentration-compactness result
\begin{Th}
\label{th-intro-3-bis}
Let $S_k\rightharpoonup S_\infty$ weakly in $\dot{W}^{1,2}(\R^2,\mbox{Sym}(n))$ where $S_k^2=id_n$ and let $u_k\rightharpoonup u_\infty$ weakly in $L^2(\R^2,{\R}^n)$  and satisfy
\[
\mbox{div}\,(S_k\,\nabla u_k)=0~~\mbox{in ${\mathcal{D}}'(\R^2)$}.
\]
Then, modulo extraction of a subsequence, there exists finitely many points $a_1\ldots a_Q\in \R^2$ s.t.
\[
u_k\longrightarrow u_\infty \quad\mbox{ strongly in }\bigcap_{p<2} W^{1,p}_{loc}(\R^2\setminus \{a_1\cdots a_Q\})\,.
\]
Moreover $u_\infty$ satisfies $\mbox{div}\,(S_\infty\,\nabla u_\infty)=0$ in ${\mathcal D}'(\R^2)$.\hfill $\Box$
\end{Th}
 We shall call $S\in \dot {W}^{1,2}(\R^2,\mbox{Sym}(n))$ where $S^2=id_n$ a {\it chirality operator}. Etymologically, in old greek $\chi\ep\iota\rho$ (kheir) means ``hand''. The word chirality refers to an intrinsic disymmetry of the space where a left and a right directions are given. More precisely almost everywhere on $\R^2$ we have the existence of two orthogonal projections, $P_R$ and $P_L$ complementary to each other ($P_R+P_L=id_n$), the left and the right,  such that $S=P_R-P_L$. 
\begin{Rm}
\label{rm-pseudo}
The system (\ref{intro-6}) is then the Euler-Lagrange equation of the Dirichlet energy into the pseudo-riemannian manifold $({\R}^n,g)$ where
\[
g(X,Y):= \lf<X,P_R Y\rg>-\lf<X,P_LY\rg>
\]
In other words (\ref{intro-6})  is the {\bf harmonic map equation} from $\R^2$ into $({\R}^n,g)$, it correspond to critical points of 
\[
E_g(u):=\int_{\R^2} |P_R\nabla u|^2-|P_L\nabla u|^2\ dx^2\,.
\]
\hfill $\Box$
\end{Rm}
As we will see theorem~\ref{th-intro-1} can be rephrased as follows.
\begin{Th}
\label{th-intro-3}
Let $P_L\in W^{1,2}(\R^2,\mbox{Sym}(n))$ such that $P_L\circ  P_L=P_L$ and denote $P_R:=id_n-P_L$ and let $f\in L^2(\R^2,{\C}^n)$ satisfying
\be
\label{intro-7}
\lf\{
\begin{array}{l}
\ds P_L\,\frac{\p f}{\p z}=0\\[5mm]
\ds P_R\,\frac{\p f}{\p \ov{z}}=0
\end{array}
\rg.
\ee
then $f\in \bigcap_{p<2}W^{1,p}_{loc}(\R^2,{\C}^n)$.\hfill $\Box$
\end{Th}
In the course of the paper we will give a third formulation of our main result. For $n=2$ it takes a simpler following form.
\begin{Th}
\label{th-intro-4}
Let $\Om\in L^2(\R^2,so(2)\otimes{\C})$ and let $f\in L^2(\R^2,{\C})$ such that
\be
\label{intro-8}
\frac{\p f}{\p z}=\Om\,\ov{f}\,.
\ee
Assume $\Im(\p_{\ov{z}}\Om)=0$, then $f\in \bigcap_{p<2}W^{1,p}_{loc}(\R^2,{\C}^2)$.\hfill $\Box$
\end{Th}
The system 
\be
\label{intro-9}
\frac{\p f}{\p z}=\Om\,{f}\,,
\ee
where $\Om\in L^2(\R^2,so(2)\otimes{\C})$ and  $\Im(\p_{\ov{z}}\Om)=0$  enjoys the same compensation property as (\ref{intro-8}) for  $f\in L^2(\R^2,{\C})$ but
this last fact is a consequence of the theory in \cite{Riv1} while theorem~\ref{th-intro-4} is new. \par
We will see that it can be recasted also in the following way.
Recall first the definition of the  {\bf Dirac Operator} in ${\C}^2$
\[
{\mathcal D}\, :=\left(
\begin{array}{cc}
0& \partial_z\\[3mm]
-\,\ov{\partial_z} &0
\end{array}
\right)
\]
Then we have the following corollary
\begin{Co}
\label{co-intro-1}
Let $U\in L^2(\R^2,{\C})$ such that $\Im(\p_{\ov{z}}U)=0$. Let $\Psi\in L^2(\R^2,{\C}^2)$ be a solution of
\be
\label{intro-10}
{\mathcal D}\Psi=\left(
\begin{array}{cc}
U& 0\\[3mm]
0 &\ov{U}
\end{array}
\right)\, \Psi\,,
\ee
then $\Psi\in \bigcap_{p<2}W^{1,p}_{loc}(\R^2,{\C}^2)$.\hfill $\Box$
\end{Co}
Throughout the paper we identify $\R^2$ with the complex number plane $\C$ and we will use both notations. \par
We will denote by ${\mathcal{S}}(\R^2)$ the space of Schwarz functions and by ${\mathcal{S}}^{\prime}(\R^2)$ the space of tempered distributions.
For $1< p < +\infty$ we will denote by $\dot{W}^{1,p}(\R^2)$ the homogeneous Sobolev space defined as the space of $f\in L^{1}_{loc}(\R^n)$ such that
$\nabla f\in L^{p}(\R^2)$ and by $\dot{W}^{-1,p^{\prime}}(\R^2)$ the corresponding dual space ($p^{\prime}$ is the conjugate of $p$).\par
We also denote by $L^{2,\infty}(\R^2)$ the space of measurable functions $f$ such that
$$
\sup_{\lambda>0}\lambda |\{x\in\R^2~: |f(x)|\ge\lambda\}|^{1/2}<+\infty\,,
$$
and $L^{2,1}(\R^2)$ is the  space of measurable functions satisfying
$$\int_{0}^{+\infty}|\{x\in\R^2~: |f(x)|\ge\lambda\}|^{1/2} d\lambda<+\infty \,\,.$$ The spaces  $L^{2,\infty}(\R^n)$ and  $L^{2,1}(\R^n)$  belongs to the family of Lorentz spaces and one  can check that they form a duality pair. For a nice introduction of   Lorentz spaces we refer to \cite{Gra}.\par
In the sequel we will often use the symbols  $a\lesssim b$ and $a\simeq b$ instead of $a\le C b$ and $C^{1} a\le b\le Cb $, whenever  the constants appearing in the estimates are not  relevant  for the computations and therefore they are omitted.

 \noindent{\bf Acknowledgments :} {\it A large part of the present work has been conceived while the two authors were visiting the Institute for Advanced Studies
in Princeton. They are very grateful to the IAS for the hospitality. The authors are also very grateful to the anonymous referee and to Jerome Wettstein for useful
remarks that permit us to improve the presentation of the paper.}

\section{Preliminaries}
  
\subsection{Bourgain-Brezis Inequalities}
In  \cite{BB}  Bourgain and Brezis proved the following striking result:
\begin{Th}[Lemma 1 in \cite{BB}]\label{ThBB}
Let $u$  be  a $2\pi$-periodic function in $\R^n$  such that $\int_{\R^n} u=0$, and let 
$\nabla u=f + g,$  where $f\in  \dot{W}^{-1,\frac{n}{n-1}}(\R^n)$ and $g\in L^1(\R^n)$  are  $2\pi$- periodic vector valued functions. Then
\begin{equation}\label{BBineq}
\|u\|_{L^{\frac{n}{n-1}}(\R^n)}\le c \left(\|f\|_{ \dot{W}^{-1,\frac{n}{n-1}}(\R^n)}+\|g\|_{L^1(\R^n)}\right).\
\end{equation}
 \end{Th}
 As a consequence of Theorem \ref{ThBB} they get the following
 \begin{Co}[Theorem 1 in  \cite{BB}]\label{corBB}
For every $2\pi$-periodic function $h\in L^n(\R^n)$  with  $\int_{R^n} h=0$ there exists a  $2\pi$-periodic $v\in W^{1,n}\cap L^{\infty}(\R^n)$ 
  satisfying
$$\div v = h~~~\mbox{in $\R^n$}$$
and
$$\|v\|_{L^{\infty}(\R^n)}+\|v\|_{W^{1,n}(\R^n)}\le C(n)\|h\|_{L^n(\R^n)}.$$
\end{Co}
\subsection{Bourgain-Brezis inequality in 2 dimension revisited}
For the convenience of the reader we provide a proof of $\eqref{BBineq}$ in $2$-dimension which has the advantage of not assuming periodicity.
 The proof is related to some compensation phenomena observed first in \cite{De} in the analysis of $2$-dimensional perfect incompressible fluids. This observation has also been used by the second author in the analysis of {\it isothermic surfaces} \cite{Riv3} (see also \cite{EM,Ma, Riv1}). 
 \par
 We start by showing the following preliminary Lemma.
 {  \begin{Lm}\label{BBLp}
Let $g\in L^1(\C,\R^2)$, $f=(f_1,f_2)\in  \dot{W}^{-1,2}(\C,\R^2)$  and  $u\in {\mathcal{S}}^{\prime}(\C,\R)$  be such that   
\begin{equation} \label{estBB-a}
\nabla u=f+g\in ( \dot{W}^{-1,2}+L^1)(\C) ~~~\mbox{in ${\mathcal{S}}^{\prime}(\C)$.}
\end{equation}
Then there is $c\in\R$ such that $u-c\in L^{2,\infty}(\C)$   and
\begin{equation}\label{L2est}
\|u-c\|_{L^{2,\infty}(\C)}\le C(\|f\|_{ \dot{W}^{-1,2}(\C)}+\|g\|_{L^1(\C)}).
\end{equation}
\end{Lm}
{\bf Proof of Lemma \ref{BBLp}.}
By assumption there exist $a_j^k\in L^2(\C,\R)$ such that
$$f_j=\sum_{k=1}^2\partial_{x_k} a^k_j.$$
For $k=1,2$ we set $a^k=(a^k_1,a^k_2)$ .   Hodge decomposition in $L^2$     gives the existence of
   $\alpha^k,\beta^k\in \dot W^{1,2}(\C)$ such that
   \begin{equation}\label{HDa}
   a^k=\nabla\alpha^k+\nabla^{\perp} \beta^k 
   \end{equation}
   and
   \begin{equation}
   \|\nabla\alpha^k\|_{L^2}+\|\nabla\beta^k\|_{L^2}\le \|a^k\|_{L^2}
   \end{equation}
\footnote{ {  One can  show \eqref{HDa}    by using the Fourier transform and the theory of Fourier symbols associated to a differential operator.
   If we denote by ${\mathcal{F}}[a^k]$ the Fourier transform of $a^k$ we have
   \begin{eqnarray*}
   {\mathcal{F}}[\partial_{x_1}(-\Delta)^{-1}(\div \,a^k)]&=&-i\xi_i|\xi|^{-2}(i\xi_1 {\mathcal{F}}[a^k_1]+i\xi_2{\mathcal{F}}[a^k_2])~~\\
  {\mathcal{F}}[-\partial_{x_2}(-\Delta)^{-1}({\rm curl} \,a^k)]&=&i\xi_2|\xi|^{-2}(-i\xi_2 {\mathcal{F}}[a^k_1]+i\xi_1{\mathcal{F}}[a^k_2])\\
  {\mathcal{F}}[\partial_{x_1}(-\Delta)^{-1}({\rm curl} \,a^k)]&=&-i\xi_1|\xi|^{-2}(-i\xi_2 {\mathcal{F}}[a^k_1]+i\xi_1{\mathcal{F}}[a^k_2]).
  \end{eqnarray*}
  Then one observes that
   \begin{eqnarray*}{\mathcal{F}}[a^k_1]&=&-i\xi_1|\xi|^{-2}(i\xi_1 {\mathcal{F}}[a^k_1]+i\xi_2{\mathcal{F}}[a^k_2])+i\xi_2|\xi|^{-2}(-i\xi_2 {\mathcal{F}}[a^k_1]+i\xi_1{\mathcal{F}}[a^k_2])\\
   &=&{\mathcal{F}}[\partial_{x_1}\alpha^k-\partial_{x_2}\beta^k]\\
   {\mathcal{F}}[a^k_2]&=&-i\xi_2|\xi|^{-2}(i\xi_1 {\mathcal{F}}[a^k_1]+i\xi_2{\mathcal{F}}[a^k_2])-i\xi_1|\xi|^{-2}(-i\xi_2 {\mathcal{F}}[a^k_1]+i\xi_1{\mathcal{F}}[a^k_2])\\
   &=&{\mathcal{F}}[\partial_{x_2}\alpha^k+\partial_{x_1}\beta^k].
   \end{eqnarray*}
   Since   $a^k\in L^2(\C)$ we have that $\nabla \alpha^k=\nabla ((-\Delta)^{-1}(\div \,a^k))\in L^{2}(\C)$, $\nabla^{\perp}\beta^k=\nabla^{\perp}((-\Delta)^{-1}({\rm curl} \,a^k))\in L^{2}(\C)$}} 
   We have
   \begin{equation} 
   \left\{\begin{array}{cc}
   \partial_{x_1} u-\sum_{k=1}^2\partial_{x_k} a^k_1=g_1&\\
    \partial_{x_2} u-\sum_{k=1}^2\partial_{x_k} a^k_2=g_2
    \end{array}\right.
    \end{equation}
    We observe that
    \begin{eqnarray}\label{dec2}
    \partial_{x_1} a^1_1&=&\partial_{x_1}(\partial_{x_1}\alpha^1-\partial_{x_2}\beta^1)\\
     \partial_{x_2} a^2_1&=&  \partial_{x_2}(\partial_{x_1}\alpha^2-\partial_{x_2}\beta^2)\\
      \partial_{x_1} a^1_2&=&\partial_{x_1}(\partial_{x_2}\alpha^1+\partial_{x_1}\beta^1) \\
       \partial_{x_2} a^2_2&=&\partial_{x_2}(\partial_{x_2}\alpha^2+\partial_{x_1}\beta^2).
       \end{eqnarray}
       Therefore we have
        \begin{equation} \label{rel}
   \left\{\begin{array}{cc}
 \partial_{x_1} u-\partial_{x_1}(\sum_{k=1}^2\partial_{x_k} \alpha^k)+\partial_{x_2}(\sum_{k=1}^2\partial_{x_k} \beta^k)=g_1& \\[5mm]
  \partial_{x_2} u-\partial_{x_2}(\sum_{k=1}^2\partial_{x_k} \alpha^k)-\partial_{x_1}(\sum_{k=1}^2\partial_{x_k} \beta^k)=g_2&  
    \end{array}\right.
    \end{equation}
       By multiplying first the second equation in \eqref{rel} and summing up the first and second one  we get
       \begin{eqnarray}\label{rel2}
      && (\partial_{x_1}+i\partial_{x_2})\left(u-\sum_{k=1}^2\partial_{x_k} \alpha^k-i(\sum_{k=1}^2\partial_{x_k} \beta^k)\right) =g_1+ig_2=:g_{\C}
      \end{eqnarray}
      By setting $w_1:=\sum_{k=1}^2\partial_{x_k} \alpha^k$ and $w_2:=\sum_{k=1}^2\partial_{x_k} \beta^k$ we have
      $$\|w_1\|_{L^2(\C)}+\|w_2\|_{L^2(\C)}\lesssim \|f\|_{\dot W^{-1,2}}.$$
and
      the equation \eqref{rel2} becomes
      
      \begin{equation}
      \partial_{\bar z} (u+w_1+iw_2)=\frac{g_{\C}}{2}
      \end{equation}
      
             We set $v:=4\pi\frac{1}{z}\ast g$. We have $\partial_{\bar z} v=g$ in ${\mathcal{S}}^{\prime}(\C).$ \footnote{We recall that $ \frac{4}{\pi z}$ satisfies $\partial_{\bar z} ( {4\pi}\frac{1}{ z})=\delta_0.$}
      Since $\frac{1}{z}\in L^{2,\infty}$ and $g\in L^1$, Young Inequality yields that $v\in L^{2,\infty}$ and
      $$\|v\|_{L^{2,\infty}}\lesssim \|\frac{1}{z}\|_{L^{2,\infty}}\|g\|_{L^1}\lesssim \|g\|_{L^1}.$$
The function
      $h=u+w_1+iw_2-v$ satisfies $\partial_{\bar z} h=0$ in ${\mathcal{S}}^{\prime}(\C)$ and therefore it is holomorphic. This implies that $\Im(h),\Re(h)$ are harmonic functions.
      By assumption $\Im(h)=w_2-\Im(v)\in L^{2,\infty}(\R^2)$ and thus $w_2-\Im(v)=0$. Since $\nabla^{\perp} (u+w_1-\Re v)=     \nabla(w_2-\Im(v))$ it follows that there is a constant $c\in\R$
      such that $u+w_1-\Re v-c=0$. This yields  in particular that $u-c\in L^{2,\infty}$.  
  The following estimate holds:
      \begin{eqnarray}
      \|u-c\|_{L^{2,\infty}}&=&\|w_1-\Re v\|_{L^{2,\infty}}\lesssim \|w_1\|_{L^{2,\infty}}+\|v\|_{L^{2,\infty}}\lesssim  \|w_1\|_{L^{2,\infty}}+C\|g\|_{L^1}\nonumber\\
      &\le&C( \|f\|_{ \dot{W}^{-1,2}(\C)}+\|g\|_{L^1(\C)}).
      \end{eqnarray}
      We can conclude the proof.~~~\hfill $\Box$\par
      \medskip
      \begin{Lm}\label{prel2}
      Let $g\in L^1(\C)$ and let $h \in L^{2,\infty}(\C)$ satisfy $\partial_{\bar z}h=g$ in ${\mathcal{S}}^{\prime}(\C)$. If  $\Im h\in L^2(\C)$, then $\Re h\in L^2(\C)$ as well and 
      \begin{equation}
      \|\Re h\|_{L^2}\le C(\|g\|_{L^1}+\|\Im h\|_{L^2}).
      \end{equation}
      \end{Lm}
      {\bf Proof of Lemma \ref{prel2}.} Let $\chi\in C^{\infty}_c(\C)$ such that $\chi=1$ on $B(0,1)$ and $\chi=0$ on $B^c(0,2).$
      For every $k\ge 1$ we set $\chi_k(x)=\chi(\frac{x}{k}).$ We set $h_k=\chi_k(x)\varphi_k\ast h$ where
      $\varphi_k\in C^\infty_c(\C)$ is a sequence of mollifiers such that $\int_{\R^2}\varphi_k dx=1.$
      We have
      \begin{equation}
      \partial_{\bar z}h_k=\chi_k\partial_{\bar z}(\varphi_k\ast h)+   \partial_{\bar z}\chi_k \varphi_k\ast h=:g_k.
      \end{equation}
      where
      $$g_k=\chi_k(\partial_{\bar z}(\varphi_k\ast h)+   \partial_{\bar z}\chi_k \varphi_k\ast h$$ and
      \begin{eqnarray}
      \|g_k\|_{L^1}&\lesssim& \|g\|_{L^1}+\frac{1}{k}[k^2]^{1/2}\|\varphi_k\ast h\|_{L^{2,\infty}}\nonumber\\
      &\lesssim& \|g\|_{L^1}+\|\varphi_k\|_{L^1}\|h\|_{L^{2,\infty}}\lesssim \|g\|_{L^1}+\|g\|_{L^1}\left \|\frac{1}{z}\right\|_{L^{2,\infty}}\nonumber\\
      &\lesssim& \|g\|_{L^1}.
      \end{eqnarray}
            For $\psi\in {\mathcal{S}}(\C)$  we define 
      \begin{eqnarray}\label{td}
      \langle \frac{1}{(\bar \xi)^2},\psi(\xi)\rangle &:=&\int_{\C}\frac{1}{(\bar \xi)^2}(\psi(\xi)-\psi(0)-\partial_{\xi_1}\psi(0)\xi_1-\partial_{\xi_2}\psi(0)\xi_2) d\xi \nonumber
\\&+&\int_{\C}\frac{1}{(\bar \xi)^2}\psi(\xi) d\xi      \end{eqnarray}
      One can see that \eqref{td} defines a tempered distribution.  
   We set $\xi=\xi_1+i\xi_2$. Observe that
 $$
 {\mathcal{F}}^{-1}\left[\frac{1}{\bar\xi^2}\right]= {\mathcal{F}}^{-1}\left[\frac{ \xi^2}{|\xi|^4}\right]=
  {\mathcal{F}}^{-1}\left[\frac{(\xi_1^2-\xi_2^2+2i \xi_1\xi_2)}{|\xi|^4}\right].
 $$
 Since ${\xi_1^2-\xi_2^2+2i \xi_1\xi_2}$ is homogeneous harmonic polynomial, we can apply Theorem 5 in 3.3 of \cite{Ste} and deduce the existence of an universal constant $c_0$  such that  
$$ {\mathcal{F}}^{-1}\lf[\frac{(\xi_1^2-\xi_2^2+2i \xi_1\xi_2)}{|\xi|^4}\rg]=c_0 \frac{x_1^2-x_2^2+ 2\, i\,x_1x_2}{|x|^2}.$$
 Now we introduce  the following tempered distribution
 
      \begin{equation}
      \hat {T}_k=-\frac{1}{(\bar \xi)^2}\hat g_k.
      \end{equation}

 We have 
   \begin{equation}
 \label{psi}
 {T}_k={\mathcal{F}}^{-1}\left[-\frac{1}{\bar\xi^2}\right]\ast g_k.
\end{equation}
  
 It follows then from (\ref{psi})
\be
\label{Linftypsi}
\|T_k\|_{L^{\infty}}\lesssim  \|g_k\|_{L^1}\ \lf\|\frac{x_1^2-x_2^2+2\, i\,x_1x_2}{|x|^2}\rg\|_{L^{\infty}}\lesssim \ \|g\|_{L^1}\,.
\ee
We also have
\begin{eqnarray}
{\mathcal{F}}[\partial_{\bar z}T_k]&=&\frac{i}{\bar \xi}{\mathcal{F}}[g_k]={\mathcal{F}}[h_k].
\end{eqnarray}
Hence
\begin{eqnarray}\label{h2}
\left|\Re\left(\int_{\C}h_k^2 dx_1dx_2\right)\right|&=&\left|\Re\left(\int_{\C}h_k\, \partial_{\bar z} T_k   dx_1dx_2\right)\right|\nonumber\\
&=& \left|-\Re\left(\int_{\C}\partial_{\bar z} h_k\, T_k   dx_1dx_2 \right)\right|= \left|-\Re\left(\int_{\C}g_k\, T_k   dx_1dx_2 \right)\right|\nonumber\\
 &\leq & \|g_k\|_{L^1}\|T_k\|_{L^{\infty}}\le C \|g\|^2_{L^1} \nonumber\\&&\end{eqnarray}
We have
\begin{equation}\label{h2h}
\Re (h^2_k)=|\Re h_k|^2-|\Im h_k|^2
\end{equation}
and 
$$\|\Im h_k\|_{L^2}^2=\int_{\C}|(\Im h \ast\varphi_k)\chi_k|^2 \lesssim\|\Im h\|_{L^2}^2.$$
From \eqref{h2} and \eqref{h2h}  we deduce that
\begin{equation}\label{ti}
\|(\Re h_k)\|^2_{L^2}\le \|\Im h_k\|_{L^2}^2+|\int_{\C}\Re h^2_k dx|\le  \|\Im h\|_{L^2}^2+C\|g\|^2_{L^1}.
\end{equation}
Up to a subsequence $\Re h_k$ converges weakly in $L^2$  to $h_{\infty}\in L^2(\R^2)$ as $k\to +\infty.$ On the other hand we have
$\Re h_k\to \Re h$ in $ {\mathcal{S}}^{\prime}(\C)$ and therefore $\Re h=h_{\infty}\in L^2(\C)$ and by the lower semicontinuity of the
$L^2$ norm we have
\begin{equation}\label{ti}
\|(\Re h)\|^2_{L^2}\le  \|\Im h\|_{L^2}^2+C\|g\|^2_{L^1}.
\end{equation}
 We conclude the proof of Lemma \ref{prel2}.~~\hfill$\Box$ 
      \par\medskip
    By combining  Lemmae  \ref{BBLp} and \ref{prel2} we can deduce  the Brezis-Bourgain Inequality.\par      
{  \begin{Lm}\label{BBL}
Let $g\in L^1(\C,\R^2)$ and   $f\in  \dot{W}^{-1,2}(\C,\R^2)$  Let $u\in {\mathcal{S}}^{\prime}(\C,\R)$  be such that   
\begin{equation} \label{estBB-b}
\nabla u=f+g\in ( \dot{W}^{-1,2}+L^1)(\C) ~~~\mbox{in ${\mathcal{S}}^{\prime}(\C)$.}
\end{equation}
Then there is $c\in\R$ such that $u-c\in L^2(\C)$   and
\begin{equation}\label{L2est}
\|u-c\|_{L^2(\C)}\le C(\|f\|_{ \dot{W}^{-1,2}(\C)}+\|g\|_{L^1(\C)}).
\end{equation}
\end{Lm}
\noindent{\bf Proof of Lemma \ref{BBL}.}

 From   Lemma  \ref{BBLp} it follows that   there is $c\in\R$ such that $u-c\in L^{2,\infty}(\C)$ with
\begin{equation}\label{L2infty}
\|u-c\|_{L^{2,\infty}(\C)}\leq C (\|f\|_{ \dot{W}^{-1,2}(\C)}+\|g\|_{L^1(\C)}).
\end{equation}

 {\bf Claim 1:}  $u-c\in L^2$ and  \eqref{L2est} holds.\par
 
 {\bf Proof of Claim 1} \par
  In the proof of Lemma \ref{BBLp} we have seen the existence of $\alpha^k,\beta^k\in\dot W^{1,2}(\C)$ ($k=1,2$) such that 
  if we set $w_1=\sum_{k=1}^2\partial_{x_k} \alpha^k$ and $w_2=\sum_{k=1}^2\partial_{x_k} \beta^k$ we have 
 $$ 
\partial_{\bar z}((u-c)-w_1-i w_2) =\frac{g_{\C}}{2} ~~~\mbox{in ${\mathcal{S}}^{\prime}(\C)$}
 $$
and 
 $$
 \|w_1\|_{L^2(\C)}+\|w_2\|_{L^2(\C)}\le C\|f\|_{ \dot{W}^{-1,2}(\C)}.
 $$
  
  The function  $h=(u-c)-w_1-i w_2$ satisfies the assumptions of Lemma \ref{prel2}. 
  Therefore we have that $(u-c)-w_1\in L^2$ with
  
 \begin{equation}\label{h3}
\int_{\C}|u-c-w_1|^2 dx \le \int_{\C} |w_2|^2dx+C\|g\|^2_{L^1} 
\end{equation}
and
\begin{eqnarray}\label{h4}
\int_{\C}|u-c|^2 dx &\lesssim &2 \|w_1\|^2_{L^2(\C)}+2 \|w_2\|^2_{L^2(\C)}+C\|g\|^2_{L^1} \nonumber\\
&\le& C(\|f\|^2_{ \dot{W}^{-1,2}(\C)}+\|g\|^2_{L^1}).
\end{eqnarray}

   We conclude the proof.~\hfill $\Box$}
 \par
 \medskip
 \begin{Rm}
 We observe that if in the  Lemma \ref{BBL} $\nabla u=\nablap v+g$ with $g\in L^1$ and $v\in L^2$ then we simply get  the estimate
  \begin{eqnarray}\label{h5}
\int_{\C}|u-c|^2 dx &\leq & \|v\|^2_{L^2(\C)}+C\|g\|^2_{L^1}
\end{eqnarray}
namely the constant in front of $\|v\|^2_{L^2(\C)}$ is $1.$
\end{Rm}

\section{Regularity of solutions to $\div(S\,\nabla u)=0$ : Proof of theorem~\ref{th-intro-1}. }
In this section we are going to investigate the regularity of $L^2$ solutions to the following system
\begin{equation}\label{maineq}
\div(S\,\nabla u)=0~~\mbox{in ${\mathcal{D}}^{\prime}(\C)$}
\end{equation}
where $S\in \dot W^{1,2}(\C,O(n))$ with  $S^2=Id$.
 \par
 It has been shown in \cite{JMS} that there exists solutions $u\in W^{1,1}_{loc}(B(0,1))$ of
 $\div(A\nabla u)=0$ in  ${\mathcal{D}}^{\prime}(B(0,1))$ where $A$ is a uniformly elliptic and continuous  matrix  which is in none of the spaces $W^{1,p}_{loc}(B(0,1))$ for any $p>1$.
 
 \medskip
 
 Actually they construct a  counter-example of a matrix $A$ which turns out to be also in $W^{1,2} (B(0,1)).$
 The matrix $A(x)=(a_{ij}(x))_{\substack{1\le i\le n\\1\le j\le n}}$  is defined as follows
 $$
 a_{ij}(x)=\delta_{ij}+\alpha(|x|)\left(\delta_{ij}-\frac{x_ix_j}{|x|^2}\right)
 $$ 
 where
 \begin{equation}\label{alpha}
 \alpha(r)=\frac{-\beta n}{(n-1)\left(\log \frac{r_0}{r}\right)}+\frac{\beta(\beta+1)}{(n-1)\left(\log\frac{r_0}{r}\right)^2}.
 \end{equation}
 where $r_0$ is large enough so that $\alpha\ge -\frac{1}{2}$ and $\beta>1$.
 
 \medskip
 
Clearly  $a_{ij}\in L^2 (B(0,1))$.  A direct computation for any $i,j,k$ gives
 
 \begin{eqnarray}\label{derivaij}
 \frac{\partial a_{ij}}{\partial {x_k}}=\alpha^{\prime}(|x|)\frac{x_k}{|x|}\left(\delta_{ij}-\frac{x_ix_j}{|x|^2}\right)
 -\alpha(|x|)\left( \frac{(\,\delta_{ik}\,x_j+\delta_{jk}\, x_i)|x|^2-2\,x_k\,x_ix_j}{|x|^4} \right).
 \end{eqnarray}
 Therefore
 $$
\left|\frac{\partial a_{ij}}{\partial{x_k}}(x)\right|\le C\frac{1}{r}\frac{1}{\log(\frac{r_0}{r})}\,.
 $$
  Since $\frac{1}{r}\frac{1}{\log(\frac{r_0}{r})}\in L^2(B(0,1))$ then $\nabla a_{ij}\in L^2 (B(0,1))$ as well.
 It is proved in \cite{JMS} that
 \begin{equation}\label{counterex}
 u(x)=x_1\frac{1}{r^2\log(\frac{r_0}{r})^{\beta}}\in L^2(B(0,1)) \quad \mbox{ solves }\quad\sum_{ij=1}^2\p_{x_i}(a_{ij}\,\p_{x_j}u)=0.\end{equation}
 The function $u$ defined in \eqref{counterex} is not in   the spaces $W^{1,p}_{loc}(B(0,1))$ for any $p>1$.
 \par
 \bigskip
  We are now proving the following result
\begin{Th}\label{S}
There is an $\varepsilon_0>0$ such that if  $S\in \dot W^{1,2}(\C, O(n))$ with $S^2=I_n$ and $\|\nabla S\|_{L^2(\C)}\le \varepsilon_0 $ then there is    $Q\in \dot{W}^{1,2}(\C,SO(n))$ such that
$$S=Q\,S^0\,Q^{-1}$$
where
\begin{equation}\label{S0g}
\ds S^0=\left(\begin{array}{c|c}
\ds I_{m\times m} &\ds 0_{m\times n-m}\\\hline 
\ds 0_{n-m\times m} &\ds -I_{n-m\times n-m}\end{array}\right)
\end{equation}
with $m\le n$
and 
\begin{equation}\label{QS}
\|\nabla Q\|_{L^2}\le C\,\|\nabla S\|_{L^2} 
\end{equation}
where $C>0$ only depends on $n$.
\hfill $\Box$
\end{Th}
\noindent{\bf Proof of Theorem \ref{S}.}  
Let $S\in \dot W^{1,2}(\R^2, O(n))$ be with $S^2=I_n$.\par
We have $\det S$, ${\mbox{Trace}}(S)\in \dot {W}^{1,2}(\R^2,\Z).$ Precisely
\begin{equation}
\det S=(-1)^{n-m},~~\mbox{and}~~{\mbox{Trace}} (S)=2m-n
\end{equation}
where $m=\# $ positive eigenvalues and $n-m=\#$ negative eigenvalues (we recall that the eigenvalues of $S$ can be either $1$ or $-1$).
 Since   $\det S,{\mbox{Trace}}(S)\in \dot{W}^{1,2}(\R^2)$ it follows that
$\det S$ and ${\mbox{Trace}}(S)$ are both constant a.e. in $\R^2.$\par
 
  We set 
  \begin{equation}\label{PLR}
  P_R:=\frac{I-S}{2}~~~\mbox{ and}~~~ P_L=\frac{I+S}{2}\end{equation}
 $P_L,P_R$ are idempotent since  $(I-S)^2=S^2-2S+I=2(I-S)$ and $(I+S)^2=2(I+S)$ and the ranks of $P_L$ and $P_R$ are constant.\par
 We can see $P_L$  (resp. $P_R$) as $\dot{W}^{1,2}$ maps with values into the Grassmanian $Gr_{m}(\R^n)$ (resp.  $Gr_{n-m}(\R^n)$) of nonoriented $m$-planes (resp. $n-m$-planes) in ${\R}^n$. 
 
 By applying Lemma 5.1.4 in H\'elein book \cite{Hel} there is an $\varepsilon_0>0$ such that one can find two  $\dot{W}^{1,2}(\R^2)$ orthonormal basis $e=e_1,\ldots,e_m$ and $f_1,\ldots,f_{n-m}$ of $Im(P_L)$ and
 $Im(P_R)$ respectively such that 
 \begin{equation}
 \|\nabla e_i\|_{L^2}\le C\|\nabla P_L\|_{L^2}~~\quad\mbox{and}\quad~~ \|\nabla f_j\|_{L^2}\le C\|\nabla P_R\|_{L^2}
 \end{equation}
 for $i=1,\ldots,m$ and $j=1,\ldots,n-m.$\par
 Let $(\epsilon_k)_{k=1,\ldots,n}$ be the canonical basis of $\R^n$. Let  $Q_L \in \dot{W}^{1,2}(\R^2,{\mathcal{M}}_{n\times n})$ and $Q_R\in \dot{W}^{1,2}(\R^2,{\mathcal{M}}_{n\times n})$ be  defined by
 \begin{equation}
 Q_L=\left(\begin{array}{c|c}
Q^L &0_{n\times n-m}\\\hline
0_{n-m\times m} & 0 \end{array}\right)
\end{equation}
and
\begin{equation}
 Q_R=\left(\begin{array}{c|c}
0&0_{n\times n-m}\\\hline
0_{n-m\times m} & Q^R \end{array}\right)
\end{equation}
with
 $Q^L\in SO(m)$, $Q^R\in SO(n-m)$ and 
 $$\{Q_L\epsilon_k,~~k=1,\ldots,m\}=\{e_j,~~j=1,\ldots m\} $$
 and $$\{Q_R\epsilon_k,~~k=m+n,\ldots,n\}=\{e_j,~~j=1,\ldots n-m\}.$$
Moreover
  $$P^L=Q_L^{-1} P_L Q_L~~\mbox{ and}~~~P^R=Q_R^{-1} P_R Q_R$$
  where
  \begin{equation}\label{P0L}
P^L=\left(\begin{array}{c|c}
I_{m\times m} &0_{m\times n-m}\\\hline
0_{n-m\times m} & 0_{n-m\times n-m}\end{array}\right)
\end{equation}
and
  \begin{equation}\label{P0L-1}
P^R=\left(\begin{array}{c|c}
0_{m\times m} &0_{m\times n-m}\\\hline
0_{n-m\times m} & I_{n-m\times n-m}\end{array}\right)
\end{equation}
  We define
  \begin{equation}
 Q=\left(\begin{array}{c|c}
Q_L &0_{n\times n-m}\\\hline
0_{n-m\times m} & Q_R \end{array}\right).
\end{equation}
 By construction we have $Q^tQ=Id,$ $S^0=Q^{-1} S Q$ and 
 \begin{equation}
 \|\nabla Q\|_{L^2}\le C\|\nabla S\|_{L^2}\,.
 \end{equation}
 This concludes the proof of Theorem \ref{S}. \hfill $\Box$
 
 Next we show how theorem~\ref{th-intro-1} implies theorem~\ref{th-intro-3}.
  More precisely we establish that \eqref{maineq} is equivalent to (\ref{intro-7}) for a suitable choice of $f$ .
 \begin{Prop}\label{equivalence}
 Let $S\in W^{1,2}(\C,O(n))$ with $S^2=I_n$ and let $u\in L^2(\C,\R^n)$ be a solution of 
\begin{equation}\label{main2}
\div(S\,\nabla u)=0~~\mbox{in ${\mathcal{D}}^\prime(\C)$}
\end{equation}
Then there exists   $v\in L^{2}(\C,\R^n)$  such that $\nablap v= S\,\nabla u$ in ${\mathcal{D}}^{\prime}(\C)$.
Moreover the function  $f=u+iv$ satisfies
  \begin{equation}\label{eqsystem}
\left\{\begin{array}{c}
\ds P_L\,\frac{\partial f}{\partial z}=0~~\mbox{in ${\mathcal{D}}^\prime(\C)$}
\\[5mm]
\ds P_R\,\frac{\partial}{\partial\bar z} f=0~~\mbox{in ${\mathcal{D}}^\prime(\C).$}
\end{array}\right.,
\end{equation}
where $P_L, P_R$ are given by \eqref{PLR}.
  \end{Prop}
\noindent{\bf Proof of Proposition \ref{equivalence}.} 
Let $v\in {\mathcal{D}}^{\prime}(\C)$ be such that $\nablap v= S\nabla u$ in  ${\mathcal{D}}^{\prime}(\C)$. It holds 
  $\nabla^{\perp}v=\nabla (Su)-\nabla S \, u\in  \dot{W}^{-1,2}+L^1.$ Lemma~\ref{BBL} gives that $v\in L^2(\C).$
We have
\begin{equation}\label{systeq}
\left\{\begin{array}{c}
S\,\partial_{x_1} u=-\partial_{x_2} v\\[5mm]
S\,\partial_{x_2} u=\partial_{x_1} v\,.
\end{array}\right.
\end{equation}
Therefore
\begin{equation}
S\,\partial_{x_1}(u+iv)=i\,\partial_{x_2}(u+iv)\,.
\end{equation}
Let us introduce $f\colon\C\to \C^n$ given by $f=u+iv.$ Obviously
$f$ satisfies:
\begin{equation}\label{systeq2}\left\{\begin{array}{c}
S\partial_{x_1} f-i\,\partial_{x_2} f=0   
\\[5mm]
S\,i\,\partial_{x_2} f-\partial_{x_1} f=0 \,.
\end{array}\right.
\end{equation}
By first subtracting and then summing the two equations in \eqref{systeq2}  we deduce that \begin{equation}
(S+I)\partial_z f=0,~~(S-I)\partial_{\bar z } f=0
\end{equation}
 Therefore $f$ satisfies (\ref{eqsystem}) and we conclude the proof.~~\hfill$\Box$

\subsection{Proof of theorem~\ref{th-intro-1} : the case $n=2$}\label{casen2}
 
  In this section we focus our attention to   the case where the function $u$ takes values in ${\R}^2$, since as we will see the formulation will become simpler and maybe more enlightening .
 
Let $Q\in \dot{W}^{1,2}(\C,SO(2))$ then a classical result by Carbou gives the existence of $\al\in \dot{W}^{1,2}({\C},{\R})$ such that
\begin{equation}\label{matrQ}
Q(x)=\left(\begin{array}{cc} \cos(\alpha(x)) & -\sin(\alpha(x))\\[5mm] \sin(\alpha(x))& \cos(\alpha(x))\end{array}\right).\end{equation}
We also set
 \begin{equation}\label{S02}
S^0=\left(\begin{array}{cc}
1 &0 \\[5mm]
0  & -1\end{array}\right)
\end{equation}
Next we re-formulate the system \eqref{intro-6} in the $n=2$ case. Precisely we have
 \begin{Prop}\label{f}
 Let $S\in \dot W^{1,2}(\C,O(2))$   with $S^2=Id$ and $\|\nabla S\|_{L^2(\C)}\le \varepsilon_0 $ (with $\varepsilon_0 >0$ as in Theorem \ref{S}). Let $Q\in SO(2)$
 as in \eqref{matrQ} such that 
 $S=Q^{-1}S^0Q$ and 
 let    $u,v$ be as in the statement of Proposition \ref{equivalence}.    Then function $f\colon\C\to(\C^2)$ \begin{equation}\label{deff}
f:=S^0Qu+iQv
\end{equation} satisfies the following equation
\be\label{maineq-ff}
\partial_z f= \left(\begin{array}{cc} 0 & 1\\ [5mm]-1& 0\end{array}\right)\partial_z\alpha\bar f\, .  
\end{equation}
 
\end{Prop}
\noindent{\bf Proof of proposition~\ref{f}}  
Let $u\in L^2(\R^2)$ be a solution of \eqref{main2} and
 $v\in L^2(\C)$ be such that 
\begin{equation}\label{veq}
\nablap v=S\,\nabla u .
\end{equation}
Since $S=Q^{-1}\,S^0\, Q$ we can write \eqref{veq} as
\begin{equation}\label{veq2}
Q\nablap v=S^0\, Q\nabla u \,.
\end{equation}

We set $f_{\Re}:=S^0\,Qu$ and $f_{\Im}:=Qv$.
From the fact that $S^0\,Q\nabla u-Q\nablap v=0$ it follows 
\begin{equation}
\nabla(f_{\Re})-\nablap(f_{\Im})=\nabla(S^0\,Q)u-\nablap Q\ v=S^0\,\nabla Q\ Q^{-1}\,S^0f_{\Re}-\nablap Q\ Q^{-1}f_{\Im}\,.
\end{equation}
Therefore
\begin{eqnarray*}
\lf\{
\begin{array}{rcl}
\ds\partial_{x_1}f_{\Re}+\partial_{x_2}f_{\Im}&=&S^0\,(\partial_{x_1}Q)Q^{-1}\,S^0f_{\Re}+\partial_{x_2}Q\,Q^{-1}\,f_{\Im}\\[5mm]
\ds\partial_{x_2}f_{\Re}-\partial_{x_1}f_{\Im}&=&S^0\,(\partial_{x_2}Q)Q^{-1}\,S^0f_{\Re}-\partial_{x_1}Q\,Q^{-1}\,f_{\Im}\,.
\end{array}
\rg.
\end{eqnarray*}
We have
\begin{eqnarray*}
\nabla Q\,Q^{-1}&=& \left(\begin{array}{cc} 0 & -1\\ 1& 0\end{array}\right)\, \nabla\al\\[5mm]
S^0\,\nabla Q\,Q^{-1}S^0&=& \left(\begin{array}{cc} 0 & 1\\ -1& 0\end{array}\right)\, \nabla\al\,.
\end{eqnarray*}
We have
\begin{equation}\label{est1}
\left\{\begin{array}{c}
\partial_{x_1}f_{\Re}+\partial_{x_2}f_{\Im}= \left(\begin{array}{cc} 0 & 1\\ -1& 0\end{array}\right)\partial_{x_1}\alpha f_{\Re}+
\left(\begin{array}{cc} 0 & -1\\ 1& 0\end{array}\right)\partial_{x_2} \alpha f_{\Im}\\[5mm]
-\partial_{x_2}f_{\Re}+\partial_{x_1}f_{\Im}=\left(\begin{array}{cc} 0 & -1\\ 1& 0\end{array}\right)\partial_{x_2} \alpha  f_{\Re}+\left(\begin{array}{cc} 0 & -1\\ 1& 0\end{array}\right)\partial_{x_1} \alpha  f_{\Im} \,.
\end{array}\right.
\end{equation}
From \eqref{est1}
it follows
\begin{eqnarray}\label{est2}
\partial_{x_1}f_{\Re}+\partial_{x_2}f_{\Im}-i(\partial_{x_2}f_{\Re}-\partial_{x_1}f_{\Im}) 
&=& \left(\begin{array}{cc} 0 & 1\\ -1& 0\end{array}\right)(\partial_{x_1} \alpha-i \partial_{x_2} \alpha)(f_{\Re}-if_{\Im})\,.\nonumber
\end{eqnarray}
Hence
\begin{equation}\label{maineq3}
\partial_z f= \left(\begin{array}{cc} 0 & 1\\ -1& 0\end{array}\right)\partial_z\alpha\bar f\,.
\end{equation}
This concludes the proof of proposition~\ref{f}. \hfill$\Box$
 
 \medskip

Now  we present the regularity of the equation \eqref{maineq-ff} and therefore of \eqref{main2}.   
We would like first to explain the reasons why the equation \eqref{maineq3}  does not fall within the classical theory of systems with a $L^2$ potential.\par

Let us   represent a function $f=u+iv$ with $u=(u_1,u_2),v=(v_1,v_2)$ as
$$
f= \left(\begin{array}{c}u_1+iv_1\\ u_2+iv_2\end{array}\right)\,.$$
We observe that
the equation  \eqref{maineq-ff}  can be written as
\be
\label{est3}
\lf\{
\begin{array}{l}
\partial_z(u_1+iv_1)= \partial_z\alpha\  (u_2-iv_2)\\[5mm]
\partial_z(  u_2+iv_2)= -\partial_z\alpha \ (u_1-iv_1)\,.
\end{array}
\rg.
\ee
 The system \eqref{est3} is of the form
 \begin{equation}\label{est4-aa}
 \left\{\begin{array}{l}
\partial_z \phi= \omega\ \bar\psi\\[5mm]
\partial_z\psi= -\,\omega\ \bar\phi\,.\end{array}\right.
\end{equation}
where $\omega=\partial_z\alpha\in L^2(\C,\C)$.
The difficulty  is that in the right hand side of \eqref{est4-aa} there are the conjugate of the unknows $(\phi,\psi).$ Suppose we would have instead
a system of the form
 \begin{equation}\label{est5-aa}
 \left\{\begin{array}{l}
\partial_z \phi= \partial_z\alpha \  \psi\\[5mm]
\partial_z\psi= -\partial_z\alpha \ \phi\,.\end{array}\right.
\end{equation}
  Then the function $\Phi:=(\phi_1,\phi_2)$ solves 
\begin{equation}
\partial_z\Phi=\Omega\ \Phi
\end{equation}
where
$$\Omega= \left(\begin{array}{cc} 0 & \partial_z\alpha\\ -\,\partial_z\alpha& 0\end{array}\right)= \p_zQ\,Q^{-1}$$
 
Hence we would deduce $\p_z(Q\,\Phi)=0$ which would imply that $\Phi\in W^{1,2}_{loc}$.   Unfortunately the multiplication of $\Phi$ solving  \eqref{est4-aa} by a matrix in $SO(2)$ does not permit to absorb the potential $\Omega$ which is the case of interest in the present work.
Therefore we have to find a different Lie group that permits us to absorb the potential.

\par
To this purpose we introduce  the algebra of  {\em Quaternions}. We 
recall standard notations regarding  this algebra    that we denote by $\K:$ 
$$\K:=\{\displaystyle a+b\   {i} +c\   {j} +d\   {k},~~(a,b,c,d)\in\R^4\}, $$
where  $i, j$ and $k$ are the fundamental quaternion units satisfying $i^2=j^2=k^2=-1$ and $ij=-ji=k$, $jk=-kj=i$ and
$ki=-ik=j$. The set $\K$  of all quaternions is a vector space over the real numbers with dimension $4$. The conjugate of $\mathfrak{q}\in \K$  is the quaternion 
 $q^{*}=a-b\ i-c\ j-d\ k.$ The reciprocal of  $\frak{q}\in \K^\ast$ is $q^{-1}=\frac{\frak{q}^*}{|\frak{q}|},$ where $|\frak{q}|=\sqrt{\frak{q}\mathfrak{q}^*}$ is the norm of $\frak{q}.$\par
  
Given $\frak{q}\in \K$, $\frak{q}=q_1+q_2i+q_3j+q_4k$ we set 
$$\Pi_i(\frak{q})=q_2 i~~\mbox{and}~~\Pi_{jk}(\frak{q})=q_3j+q_4 k\,.$$
 We also denote  by ${\K}_p$ the quaternion of the form
$\frak{q}=q_2i+q_3j+q_4k$ (the {\em pure quaternions}) and ${\mathcal{U}}(\K):=\{\frak{q}\in  \K: ~~|\frak{q}|=1\}.$
  ${\K}_p$ is the Lie Algebra of the Lie Group ${\mathcal{U}}({\K}).$

Finally given $\frak{f}\colon\C\to \K$ we introduce the following differential operators  (Cauchy-Riemann-Fueter operators):
\begin{eqnarray}
\partial_L \frak{f}&:=&2^{-1}(\partial_{x_1}\frak{f}-i\ \partial_{x_2}\frak{f})\label{deltal}\\[5mm]
\partial_R f&:=&2^{-1}(\partial_{x_1}\frak{f}-\partial_{x_2}\frak{f} \ i) \label{deltaR}.
\end{eqnarray}
and
\begin{eqnarray}
\bar\partial_L \frak{f}&:=&2^{-1}(\partial_{x_1}\frak{f}+i\ \partial_{x_2}\frak{f})\label{deltal-aa}\\[5mm]
\bar\partial_R f&:=&2^{-1}(\partial_{x_1}\frak{f}+\partial_{x_2}\frak{f} \ i )\label{deltaR-aa}.
\end{eqnarray}
We observe that if $\frak{f}$ takes values in $\C$ then
$$\partial_L \frak{f}=\partial_R \frak{f}=\partial_z \frak{f}~~\mbox{and}~~\bar\partial_L \frak{f}=\bar\partial_R \frak{f}=\partial_{\bar z} \frak{f}.$$
We are going to {rewrite the equation  \eqref{maineq-ff} and therefore \eqref{est3} using the quaternion valued functions}.
\begin{Lm}
\label{lm-quaternions}
Let
$$
f= \left(\begin{array}{c}u_1+iv_1\\[5mm] u_2+iv_2\end{array}\right)$$
be a solution of \eqref{maineq-ff} then 
 the quaternion 
$$
\frak{f}=u_1+v_1i+ u_2j+v_2k$$ 
satisfies
\begin{eqnarray}\label{maineq4}
\partial_L\frak{f}&=&-\,\partial_{ z}\alpha\,j\  \frak{f}\,.
\end{eqnarray}
\end{Lm}
\noindent{\bf Proof of lemma~\ref{lm-quaternions}.}
We have seen that the equation \eqref{maineq-ff} is equivalent to the system \eqref{est3}. Such a system can also be written using the  $\partial_L$ operator, which coincides with $\partial_R$ at this stage since the variables $u_1+iu_2$ and $v_1+iv_2$ are $\C$-valued.
\be
\label{est3-bis}
\lf\{
\begin{array}{l}
\partial_L( u_1+iv_1)=\partial_z\alpha \ (u_2-iv_2)\\[5mm]
\partial_L( u_2+iv_2)=-\partial_z\alpha \ (u_1-iv_1)\,.
\end{array}
\rg.
\ee
We multiply from the right the second equation in \eqref{est3-bis} by $j$ and we get (recall that $ij=k=-ji$)
\begin{eqnarray}\label{est4}
\partial_L( u_2j+v_2k)&=&-\partial_z\alpha j (u_1+iv_1)\,.
\end{eqnarray}
On another hand we can write the   first equation in   \eqref{est3-bis} as follows:
\be\label{est5}
\partial_L( u_1+iv_1)= -\,\partial_z\alpha \, j^2\ (u_2-iv_2)=-\,\partial_z\alpha j \,(u_2\,j+v_2\,k).
\end{equation}
By summing \eqref{est4} and \eqref{est5}
we find
\begin{eqnarray}\label{est6-aa}
\partial_L( u_1+v_1i+ u_2j+v_2 k)&=&\,-\partial_z\alpha j\ (u_1+v_1i+ u_2j+v_2 k)\,.
\end{eqnarray}
Hence we get \eqref{maineq4} and we can conclude.~~\hfill$\Box$
\par
\bigskip

\subsection{Bootstrap test for {$\mathbf \partial_L\frak{f}=\partial_{z}\alpha j\ \frak{f}$}}\label{boot}
In the sequel up to exchange $\alpha$ and $-\alpha$ we   study the following equation
\begin{equation}\label{maineq5}
 \partial_L\frak{f}=\partial_{z}\alpha j\ \frak{f}
\end{equation}
Actually all that is proved in this section also holds for a system of the form
\begin{equation}\label{maineq5bis}
 \partial_L\frak{f}=\Omega j\ \frak{f}
\end{equation}
where $\Omega\in L^2(\C,\C).$\par

The first main goal of this section is to show that the operator

$$\frak{f} \in L^2(\C,\K)\mapsto \partial_L\frak{f}-\partial_{z}\alpha j\ \frak{f}$$
is {\bf injective} if the $L^2$ norm of $\partial_{z}\alpha $ is sufficiently small. This is what we call the {\bf``bootstrap test''}.
\begin{Th}\label{mainth} There exists $\varepsilon_0>0$ such that for every $\alpha\in \dot W^{1,2}(\C,\R)$ satisfying $\|\nabla\alpha\|_{L^2}\le \varepsilon_0$ and every  $\frak{f} \in L^2(\C,\K)$ solving
 \begin{equation}\label{maineq6}
\partial_L \frak{f}=\partial_{z}\alpha \ j\ \frak{f}\ ,
\end{equation}
then ${\frak{f}}\equiv 0.$ 
\end{Th}
%\noindent{\bf Proof of Theorem \ref{mainth}.}
\noindent Before going to the proof of theorem~\ref{mainth} we will introduce a nonlinear operator ${\mathbf N}$.

\medskip

Let $\frak{q}\in {\mathcal{U}}({\K})$.  We multiply the equation \eqref{maineq6} on the \underbar{left} by $\frak{q}$ :
 \begin{equation}\label{est6}
\frak{q}[\partial_{x_1} {\frak{f}}-i\ \partial _{x_2} {\frak{f}}]=\frak{q}[\partial_{x_1} \alpha-\partial _{x_2} \alpha\ i] j {\frak{f}}.
\end{equation}
Observe that
\begin{eqnarray}\label{est7}
\frak{q}[\partial_{x_1}  {\frak{f}}-\ i \partial _{x_2} {\frak{f}} ]&=&\partial_{x_1}[\frak{q} {\frak{f}}]-\partial_{x_2}[\frak{q}\ i {\frak{f}}]
-\partial_{x_1}\frak{q}\ {\frak{f}}+\partial_{x_2}\frak{q} \ i \ {\frak{f}}.
\end{eqnarray}
By combining \eqref{est6} and \eqref{est7} we get
\begin{eqnarray}\label{est8}
\partial_{x_1}[\frak{q}{\frak{f}}]-\partial_{x_2}[\frak{q} \ i {\frak{f}} ]&=& \frak{q}[\partial_{x_1} \alpha-\partial _{x_2} \alpha\ i] j {\frak{f}}\nonumber\\[5mm]
&+&\partial_{x_1}\frak{q}\ {\frak{f}}-\partial_{x_2}\frak{q}\ i\ {\frak{f}} \\[5mm]
&=&\frak{q}[\partial_{x_1} \alpha j-\partial _{x_2} \alpha k+\frak{q}^{-1}\partial_{x_1}\frak{q}-\frak{q}^{-1}\partial_{x_2}\frak{q}\ i] {\frak{f}}.\nonumber
\end{eqnarray}
\par
\medskip
We observe that since $|\frak{q}|=1$ then $\frak{q}^{-1}\partial_{x_i}\frak {q}\in \K_p,$.   \par
We introduce the following operator
\begin{eqnarray}\label{est9}
&& {\mathbf N}\colon \dot W^{1,2}(\C, {\mathcal{U}}({\K}))\to  \dot{W}^{-1,2}(\C,\mbox{Span$\{i\}$})\times L^{2}(\C,\mbox{Span$\{j,k\}$})\nonumber\\[5mm]
&&~~~~~ \frak{q}\mapsto(\Pi_i(\partial_{x_1}(\frak{q}^{-1}\partial_{x_1} \frak{q})+\partial_{x_2}(\frak{q}^{-1}\partial_{x_2} \frak{q})),\Pi_{jk}(\frak{q}^{-1}\partial_{x_1} \frak{q}-\frak{q}^{-1}\partial_{x_2} \frak{q}\,i)
\end{eqnarray}
 We shall prove the following result.
\begin{Lm}\label{lm-invertN-a}
There is  $\varepsilon_0>0$ and $C>0$ such that for any choice of $\om\in  \dot{W}^{-1,2}({\C},{i\R})$ and $\frak{g}\in L^2({\C},\mbox{Span}\{j,k\})$ satisfying
\begin{equation}\label{smallass}\|\omega\|_{ \dot{W}^{-1,2}}\le \varepsilon_0,~~\|\frak{g}\|_{L^{2}}\le \varepsilon_0
\end{equation}
then
there is $\frak{q}\in \dot  W^{1,2}(\C,{\mathcal{U}}({\K}))$.
\begin{equation}\label{decomp}{\mathbf N}(\frak{q})=(\omega,\frak{g})\end{equation}
and
\be
\|\nabla\frak{q}\|_{L^2}\le C(\|\omega\|_{ \dot{W}^{-1,2}}+\|\frak{g}\|_{L^{2}})\,.
\ee
\end{Lm}
In order to prove lemma~\ref{lm-invertN-a} we shall need to introduce some notations and establish some intermediate results.\par
As in \cite[Proof of Theorem~1.2, Step 4]{DLR2}, by an approximation argument it suffices to prove Lemma~\ref{lm-invertN-a} assuming that $\om$ and $\frak{g}$ are slightly more regular. \par 
We fix $2<p<+\infty$ and for $\varepsilon > 0$ we introduce
\be
\label{Uepsilon}
     \mathcal{U}_\eps  := \left \{
     \begin{array}{c}
     (\omega,\frak{g}) \in  \dot{W}^{-1,p}\cap  \dot{W}^{-1,p^{\prime}}(\C,i\R)\times L^p\cap L^{p^{\prime}}(\C,\mbox{Span$\{j,k\}$})\\[5mm]
   \|\omega\|_{ \dot{W}^{-1,2}} +\|\frak{g}\|_{L^2}\le \varepsilon
  \end{array}
  \right \}
  \ee
where $\ds\frac{1}{p}+\frac{1}{p'}=1$. \par
For constants $\varepsilon, \Theta > 0$ let $\mathcal{V}_{\varepsilon,\Theta} \subseteq \mathcal{U}_\varepsilon$ be the set where we have the decomposition \eqref{decomp} with the estimates
 \be
\label{eq:gauge:2est}
  \|\nabla \frak{q}\|_{2} \leq \Theta  ( \|\omega\|_{ \dot{W}^{-1,2}}+\|\frak{g}\|_{L^2})\, 
 \ee 
 \be
  \label{eq:gauge:pest}
 \|\nabla \frak{q}\|_{p} \leq \Theta  ( \|\omega\|_{ \dot{W}^{-1,p}}+\|\frak{g}\|_{L^p})\, ,
 \ee
 \be
 \label{eq:gauge:ppest}
     \|\nabla \frak{q}\|_{p^{\prime}} \leq \Theta  ( \|\omega\|_{ \dot{W}^{-1,p^{\prime}}}+\|\frak{g}\|_{L^{p^\prime}})\,.
\ee
That is
\[
 \mathcal{V}_{\varepsilon,\Theta} := \left \{\omega,\frak{g} \in \mathcal{U}_\varepsilon:\ \begin{array}{c}
                                             \mbox{there exists $\frak{q} \in (\dot{W}^{1,p} \cap \dot{W}^{1,p^{\prime}})(\R^2,{\mathcal{U}}({\K}))$, so that $\frak{q}-\frak{1}\in L^{\frac{2p}{p-2}}$}\\[3mm]
                                            \mbox{and}~~ \eqref{decomp}, \eqref{eq:gauge:2est}, \eqref{eq:gauge:pest},  \eqref{eq:gauge:ppest}~~ \mbox{hold.}
                                              \end{array}\right \}
\]
 \footnote{Note that \eqref{eq:gauge:2est} could actually be deduced  from \eqref{eq:gauge:pest},  \eqref{eq:gauge:ppest} by interpolation.}
The strategy   to prove lemma~\ref{lm-invertN-a} follows the one K. Uhlenbeck introduced in \cite{Uh} to construct Coulomb gauges in critical dimensions. 
Precisely  lemma~\ref{lm-invertN-a}  is going to be a consequence of the following proposition.
\begin{Prop}\label{pr:VepseqUeps}
There exist $\Theta > 0$ and $\eps > 0$ so that $\mathcal{V}_{\eps,\Theta} = \mathcal{U}_\eps$. \hfill $\Box$
\end{Prop}
\noindent{\bf Proof of Proposition \ref{pr:VepseqUeps}.}
Proposition~\ref{pr:VepseqUeps} follows, once we show the following four properties
\begin{itemize}
 \item[(i)] $\mathcal{U}_\eps$ is  connected.
 \item[(ii)] $\mathcal{V}_{\eps,\Theta}$ is nonempty. 
 \item[(iii)] For any $\eps, \Theta > 0$, $\mathcal{V}_{\eps,\Theta}$ is a relatively closed subset of $\mathcal{U}_\eps$.
 \item[(iv)] There exist $\Theta > 0$ and $\eps > 0$ so that $\mathcal{V}_{\eps,\Theta}$ is a relatively open subset of $\mathcal{U}_\eps$.
\end{itemize}
\par
Property (i) is clear, since $\mathcal{U}_\eps$ is obviously starshaped with center $0$.

\medskip

Property (ii) is also obvious since : $(0,0)\in \mathcal{V}_{\eps,\Theta}$. 

\medskip

The closedness property (iii) follows almost verbatim as in  \cite[Proof of Theorem~1.2, Step~1, p.1315]{DLR2}: there one replaces $(-\Delta)^{1/4}$ by $\nabla$.
  Observe that a uniform bound of the $L^{p}$-norm as in \eqref{eq:gauge:ppest} implies by Sobolev embedding in particular a uniform bound of $\frak{q}-\frak{1}$ in $L^{\frac{2p}{p-2}}(\R^2)$.
  
\medskip

We show now the {openness property} (iv). For this let $\omega_0,\frak{g}_0 $ be arbitrary in $\mathcal{V}_{\eps,\Theta}$, for some $\eps, \Theta > 0$ chosen below. \par Let $\frak{q}_0 \in \dot{W}^{1,p} \cap \dot{W}^{1,p^{\prime}}(\C,{\mathcal{U}}({\K}))$, $\frak{q}_0-1 \in L^{\frac{2p}{p-2}}(\C )$ so that the decomposition      \eqref{decomp} as well as the estimates  \eqref{eq:gauge:2est}, \eqref{eq:gauge:pest} and \eqref{eq:gauge:ppest}  are satisfied for $\omega_0$ and $\frak{g}_0 $.

\medskip

We consider perturbations of $\frak{q}_0$ of the form $\frak{q}=\frak{q}_{0}e^{\frak{u}}$  
where $\frak{u}\in (\dot W^{1,p} \cap \dot W^{1,p^{\prime}} \cap L^{\frac{2p}{p-2}})(\C,{\K}_p).$ Observe that the exponent $p>2$ has been chosen in particular to ensure $\frak{u}\in (C^0\cap L^\infty)({\C})$ with uniform estimates and $\frak{q}_{0}e^{\frak{u}}-\frak{1}\in L^{\frac{2p}{p-2}}$.\vfill\eject
\footnote{Indeed for a Schwartz function one has
\[
\frak{u}(x)=\frac{1}{2\pi}\int_{\C}\nabla_x\log|x-y|\cdot\nabla\frak{u}(y)\ dy\quad\Rightarrow\quad\|\frak{u}\|_\infty\le (2\pi)^{-1}\|\nabla_x\log|x-y|\|_{L^{2,\infty}}\, \|\nabla\frak{u}\|_{L^{2,1}}
\] 
Generalized H\"older inequality (see \cite{Gra}) gives moreover
\[
 \|\nabla\frak{u}\|_{L^{2,1}}\le C\,  \|\nabla\frak{u}\|^{\alpha}_{L^p}\  \|\nabla\frak{u}\|^{1-\alpha}_{L^{p^\prime}}\,.
\] 
where $2^{-1}=\alpha p^{-1}+(1-\alpha){p'}^{-1}$. If $\frak{u}\in  \dot W^{1,p} \cap \dot W^{1,p^{\prime}} \cap L^{\frac{2p}{p-2}}$ then
$\tilde{\frak{u}}(x)=\frac{1}{2\pi}\int_{\C}\nabla_x\log|x-y|\cdot\nabla\frak{u}(y)\ dy$ satisfies $\Delta \tilde{\frak{u}}=\Delta \frak{u}$ in ${\mathcal{S}}^{\prime}(\C)$. 
One has 
$\|\tilde{\frak{u}}\|_\infty\le (2\pi)^{-1}\|\nabla_x\log|x-y|\|_{L^{2,\infty}}\, \|\nabla\frak{u}\|_{L^{2,1}}$
 and moreover $\tilde{\frak{u}}\in L^{\frac{2p}{p-2}}.$   $\tilde{\frak{u}}-\frak{u}$ is a harmonic function belonging to $L^{\frac{2p}{p-2}}$, hence $\tilde{\frak{u}}-\frak{u}=0.$}
 \par
 We set  
\begin{eqnarray}\label{tilden}
\tilde{N}_{\frak{q}_{0}}\colon  (\dot W^{1,p} \cap \dot W^{1,p^{\prime}} \cap L^{\frac{2p}{p-2}})(\C,{\K}_p)&\to& \left(i( \dot{W}^{-1,p}\cap  \dot{W}^{-1,p^{\prime}})(\C), L^p\cap L^{p^\prime}(\C,\mbox{span$\{j,k\}$})\right)\nonumber
\\
 \frak{u}&\mapsto&\tilde{N}_{\frak{q}_{0}}(\frak{u}):= {\mathbf{N}}(\frak{q}_{0}e^{\frak{u}}). 
\end{eqnarray}
We will write  $\frak{u}=\frak{u}_1i+\frak{u}_2j+\frak{u}_3 k.$

\medskip

We have $\tilde{N}_{\frak{q}_{0}}\in C^1$ and   we can compute $D \tilde{N}_{\frak{q}_{0}}(0)$ as 
\[
D\tilde{N}_{\frak{q}_{0}}(0)= \frac{d}{dt} \tilde{N}_{\frak{q}_{0}}(t{\frak{u}}) \Big|_{t=0}=: L_{\frak{q}_{0}}({\frak{u}}),
\]
where for ${\frak{u}} \in L^{\frac{2p}{p-2}}\cap \dot W^{1,p} \cap \dot W^{1,p^{\prime}}(\C,{\K}_p)$ 
\be
\label{diff}
\begin{array}{l}
\ds L_{\frak{q}_{0}}({\frak{u}}) :=\left(\Pi_i\left(\Delta {\frak{u}}+\partial_{x_1}[\frak{q}_{0}^{-1}\partial_{x_1}\frak{q}_{0} {\frak{u}}-{\frak{u}} \frak{q}_{0}^{-1}\partial_{x_1}\frak{q}_{0}]+\partial_{x_2}[\frak{q}_{0}^{-1}\partial_{x_2}\frak{q}_{0} {\frak{u}}-{\frak{u}} \frak{q}_{0}^{-1}\partial_{x_2}\frak{q}_{0}]\right),\right.\\[5mm]
\ds\quad\quad\quad\left.\Pi_{jk}(\partial_{x_1}{\frak{u}}-\partial_{x_2}{\frak{u}} i+[\frak{q}_{0}^{-1}\partial_{x_1}\frak{q}_{0} {\frak{u}}-{\frak{u}} \frak{q}_{0}^{-1}\partial_{x_1}\frak{q}_{0}]-[\frak{q}_{0}^{-1}\partial_{x_2}\frak{q}_{0} {\frak{u}}-{\frak{u}} \frak{q}_{0}^{-1}\partial_{x_2}\frak{q}_{0}]i)\right)\,.
\end{array}
 \ee
In order to use a fixed-point argument for $\tilde N_{\frak{q}_{0}}$, we will show that $L_{\frak{q}_{0}}$ is an isomorphism. More precisely we prove the following lemma.

\begin{Lm}\label{la:invertible}
For any $\Theta > 0$ there exists a $\eps > 0$ so that the following holds for any $\omega_0,\frak{g}_0 $  and $\frak{q}_0$ as above.
\par
For any $\omega \in ( \dot{W}^{-1,p}\cap  \dot{W}^{-1,p^{\prime}})(\C,i\R)$, $\frak{g}\in ( L^p\cap L^{p^\prime})(\C,\mbox{span$\{j,k\}$})$ there exists a unique ${\frak{u}}\in L^{\frac{2p}{p-2}}\cap \dot W^{1,p} \cap \dot W^{1,p^{\prime}}(\C,{\K}_p)$  so that 
\[
(\omega,\frak{g}) = L_{\frak{q}_0}({\frak{u}}) 
\]
and for some constant $C = C(\omega_0,\alpha_0,\Theta) > 0$ it holds
\begin{eqnarray}\label{esteta}
\|{\frak{u}}\|_{L^{\frac{2p}{p-2}}} + \|\nabla{\frak{u}}\|_{L^{p}(\C)} +\|\nabla{\frak{u}}\|_{L^{p^{\prime}}(\C)} &\le&
 C\left[ \brac{\|\omega\|_{ \dot{W}^{-1,p}(\C)} +\|\omega\|_{ \dot{W}^{-1,p^{\prime}}(\C)}}\right.\nonumber\\[3mm]
 &+&\left.\brac{\|\frak{g}\|_{L^{p}(\C)} +\|\frak{g}\|_{L^{p^{\prime}}(\C)}}\right].
\end{eqnarray}
\end{Lm}
\noindent{\bf Proof of lemma~\ref{la:invertible}.}

\noindent{\bf Claim 1.}  ${\mathbf L_{1}({\frak{u}})}$ {\bf is invertible} $(\frak{q}_0=1)$ as a map    
$$
L_1\colon  ( \dot {W}^{1,p^{\prime}} \cap L^{\frac{2p}{p-2}})(\C,{\K}_p)
\rightarrow \left(  \dot{W}^{-1,p^{\prime}}(i\R)\times (L^p\cap L^{p^\prime})(\C,\mbox{span$\{j,k\}$})\right)$$
The operator $L_1$ is given by
\be
\label{est17}
\begin{array}{l}
\ds L_{1}({\frak{u}})=\frac{d}{dt}\mathbf{N}(e^{t\frak{u}})_{t=0}=\left(\Pi_{i}(\Delta \frak{u}),\Pi_{jk}\left(\partial_{x_1} \frak{u}-\partial_{x_2}\frak{u} i\right)\right)\\[5mm]
\ds\quad\quad\quad=\left(\Delta \frak{u}_1\ i,(\partial_{x_1} \frak{u}_2-\partial_{x_2}\frak{u}_3)j+(\partial_{x_1} \frak{u}_3+\partial_{x_2}\frak{u}_2)k\right)\,.
\end{array}
\ee
  Given $f\in  \dot{W}^{-1,p'}(\C,\R)$, $a,b\in  L^{p^\prime}(\C,\R)$ there is a unique triple $\frak{u}_1,\frak{u}_2,\frak{u}_3\in  \dot W^{1,p} \cap\dot W^{1,{p^\prime}} \cap L^{\frac{2p}{p-2}}({\C},\R) $ such that
$$L_{1}({\frak{u}})=(fi ,aj+bk)\,.$$
More precisely the following system should be satisfied:
\begin{equation}\label{est18}
\left\{\begin{array}{l}
\Delta \frak{u}_1=f\\[5mm]
\partial_z(\frak{u}_2-\frak{u}_3\,i)=a-bi\end{array}\right.
\end{equation}
$$
\frak{u}_1(x)=\frac{1}{2\pi}\log(|x|)\ast f(x)\,,\quad\frak{u}_2-\frak{u}_3i=\frac{1}{4\pi}(a-bi)\ast \frac{1}{\bar z}\,.$$
Classical estimates give
$$\|\frak{u}_1\|_{L^{\frac{2p}{p-2}}}+\|\nabla \frak{u}_1\|_{L^{p^\prime}}\lesssim\|f\|_{  \dot{W}^{-1,p^{\prime}} },~~\|\frak{u}_2-\frak{u}_3\,i\|_{L^{\frac{2p}{p-2}}}+\|\nabla(\frak{u}_2-\frak{u}_3\,i)\|_{L^{p^\prime}}\lesssim \|a-b\,i\|_{L^{p^\prime}}\,.$$
{\bf The Claim 1 is proved.}

\medskip

Observe that
\[
\begin{array}{l}
\ds L_{\frak{q}_{0}}({\frak{u}})-L_1({\frak{u}})=\left(\Pi_{i}(\partial_{x_1}[\frak{q}_{0}^{-1}\partial_{x_1}\frak{q}_{0} {\frak{u}}-{\frak{u}} \frak{q}_{0}^{-1}\partial_{x_1}\frak{q}_{0}]+\partial_{x_2}[\frak{q}_{0}^{-1}\partial_{x_2}\frak{q}_{0} {\frak{u}}-{\frak{u}} \frak{q}_{0}^{-1}\partial_{x_2}\frak{q}_{0}]),\right.\\[5mm]
\ds\quad\quad\quad\left.\Pi_{jk}([\frak{q}_{0}^{-1}\partial_{x_1}\frak{q}_{0} {\frak{u}}-{\frak{u}} \frak{q}_{0}^{-1}\partial_{x_1}\frak{q}_{0}]-[\frak{q}_{0}^{-1}\partial_{x_2}\frak{q}_{0} {\frak{u}}-{\frak{u}} \frak{q}_{0}^{-1}\partial_{x_2}\frak{q}_{0}]i)\right)\,.
\end{array}
\]
We have
\be\label{esteta-aa}
\begin{array}{l}
\ds\|\partial_{x_\ell}\left( \frak{q}_{0}^{-1}\partial_{x_\ell}\frak{q}_{0}\,{\frak{u}}-{\frak{u}}\, \frak{q}_{0}^{-1}\partial_{x_\ell}\frak{q}_{0}\right)\|_{ \dot{W}^{-1,p^\prime}}\le \|\frak{q}_{0}^{-1}\partial_{x_\ell}\frak{q}_{0}\,{\frak{u}}-{\frak{u}} \,\frak{q}_{0}^{-1}\partial_{x_\ell}\frak{q}_{0}\|_{L^{p^\prime}}\\[5mm]
\quad\le  \|\nabla \frak{q}_{0}\|_{L^2}\, \|{\frak{u}}\|_{L^{\frac{2p}{p-2}}}\le \varepsilon\ \Theta\ \|{\frak{u}}\|_{L^{\frac{2p}{p-2}}}
\end{array}
\ee
Choosing $\varepsilon>0$ small enough (depending on $\Theta$) we obtain that $L_{\frak{q}_{0}}$ is an invertible map from $ \dot W^{1,p^\prime}\cap  L^{\frac{2p}{p-2}}$ to $ \dot{W}^{-1,p^\prime}(\C,i\R)\times L^{p'} (\C,\mbox{span$\{j,k\}$}),$ 

\medskip
\noindent{\bf Claim 2.} Assuming now $\omega \in ( \dot{W}^{-1,p}\cap  \dot{W}^{-1,p^{\prime}})(\C,i\R)$, $\frak{g}\in ( L^p\cap L^{p^\prime})(\C,\mbox{span$\{j,k\}$})$ we prove that
the unique solution  ${\frak{u}}$ of  $L_{\frak{q}_{0}}({\frak{u}})$ is in  $\dot W^{1,p}$.

 From the fact that $(w,\frak{g})=L_{\frak{q}_{0}}$ it follows:  \par
 \begin{equation}\label{estgauge}
 \left\{\begin{array}{c}
 {\Delta{\frak{u}_1}}i= \omega -\,\Pi_{i}(\partial_{x_1}\left( \frak{q}_{0}^{-1}\partial_{x_1}\frak{q}_{0}{\frak{u}}-{\frak{u}} \frak{q}_{0}^{-1}\partial_{x_1}\frak{q}_{0}\right))
 -\Pi_{i}(\partial_{x_2}\left( \frak{q}_{0}^{-1}\partial_{x_2}\frak{q}_{0}{\frak{u}}-{\frak{u}} \frak{q}_{0}^{-1}\partial_{x_2}\frak{q}_{0}\right))
 \\[5mm]
(\partial_{x_1} \frak{u}_2-\partial_{x_2}\frak{u}_3)j+(\partial_{x_1} \frak{u}_3+\partial_{x_2}\frak{u}_2)k=\frak{g}\\[5mm]~~~~~~+\,\Pi_{jk}\left(-[\frak{q}_{0}^{-1}\partial_{x_1}\frak{q}_{0} {\frak{u}}-{\frak{u}} \frak{q}_{0}^{-1}\partial_{x_1}\frak{q}_{0}]+[\frak{q}_{0}^{-1}\partial_{x_2}\frak{q}_{0} {\frak{u}}-{\frak{u}} \frak{q}_{0}^{-1}\partial_{x_2}\frak{q}_{0}]i\right).\end{array}\right.
 \end{equation}
 We observe that  
 \begin{eqnarray*}
 && (\partial_{x_1} \frak{u}_2-\partial_{x_2}\frak{u}_3)j+(\partial_{x_1} \frak{u}_3+\partial_{x_2}\frak{u}_2)k= {2\partial_{\bar z}(\frak{u}_2+ \frak{u}_3\,i)j}.\end{eqnarray*}
 Therefore we can write the second equation in \eqref{estgauge} in the following way:
 \begin{eqnarray}\label{estgaugebis}
2 \partial_{\bar z}( \frak{u}_2+ \frak{u}_3\,i)&=&- \frak{g}j +\Pi_{jk}[\frak{q}_{0}^{-1}\partial_{x_1}\frak{q}_{0} {\frak{u}}-{\frak{u}} \frak{q}_{0}^{-1}\partial_{x_1}\frak{q}_{0}]j+\left(\Pi_{jk}([\frak{q}_{0}^{-1}\partial_{x_2}\frak{q}_{0} {\frak{u}}-{\frak{u}} \frak{q}_{0}^{-1}\partial_{x_2}\frak{q}_{0}]i\right)j.\nonumber\\&&
 \end{eqnarray}
 
 Let $p'<r<2$, since $\nabla \frak{q}_0\in L^p$ we have for $\ell=1,2$
 \begin{equation}\label{estgauge2}
 \|\frak{q}_{0}^{-1}\partial_{x_\ell}\frak{q}_{0}{\frak{u}}\|_{L^t}\lesssim \|\nabla \frak{q}_0\|_{L^p}\|{\frak{u}}\|_{L^{\frac{2r}{2-r}}}
 \end{equation}
 for $\frac{1}{t}=\frac{1}{p}+\frac{2-r}{2r}.$ Observe that $p>t>2,$ since $r>p^{\prime}.$\par
 From \eqref{estgauge} and \eqref{estgauge2} it follows $\nabla {\frak{u}}\in L^t$. We have also $\nabla  {\frak{u}}\in L^{t^\prime}$. This  implies that ${\frak{u}}\in L^{\infty}$ (see previous footnote). Therefore
 \begin{equation}\label{estgauge3}
 \|\frak{q}_{0}^{-1}\partial_{x_\ell}\frak{q}_{0}{\frak{u}}\|_{L^p}\lesssim \|\nabla \frak{q}_0\|_{L^p}\|{\frak{u}}\|_{\infty}
 \end{equation}
 From \eqref{estgauge} it follows that $\nabla {\frak{u}}\in L^{p}$  and the Claim 2 is  proved. This concludes the proof of lemma~\ref{la:invertible}.
\hfill $\Box$

 \medskip
 
 \noindent{\bf Proof of   proposition \ref{pr:VepseqUeps} continued.}

  \par
 For $\eps = \eps(\Theta) > 0$ chosen small enough, and for any $(\omega_0,\frak{g}_0) \in \mathcal{V}_{\eps,\Theta}$ the local inversion theorem applied to ${\mathbf{N}}$ gives the existence of some $\delta > 0$ (that might depend on $(\omega_0,\frak{g}_0)$) such that, for every $(\omega,\frak{g}) \in \mathcal{U}_\eps$ with

\begin{eqnarray}
&& \|\omega - \omega_0 \|_{ \dot{W}^{-1,p}(\C)}+\|\omega - \omega_0 \|_{ \dot{W}^{-1,p^{\prime}}(\C)} < \delta\label{estgauge5}\\[5mm]
&& \|\frak{g} - \frak{g}_0 \|_{L^{p}(\C)}+\|\frak{g}- \frak{g}_0 \|_{L^{p^{\prime}}(\C)} < \delta, \label{estgauge6}
\end{eqnarray}
we find $\frak{q} = \frak{q}_0 e^{\frak{u}} \in \dot{W}^{1,p} \cap \dot{W}^{1,p^{\prime}}(\C,{\mathcal{U}}({\K}))$, so that $\frak{q}-1 \in L^{\frac{2p}{p-2}}(\C,{\K})$ and \eqref{decomp} is satisfied.

It remains to prove \eqref{eq:gauge:2est}, \eqref{eq:gauge:pest} and \eqref{eq:gauge:ppest}. The local inversion theorem does not imply the estimates
\eqref{eq:gauge:2est}, \eqref{eq:gauge:pest}
                                            and \eqref{eq:gauge:ppest}. Anyway for every $\omega,\frak{g}\in {\mathcal{U}}_{\varepsilon}$  satisfying  \eqref{estgauge5} and \eqref{estgauge6} and for every $\sigma$ we can choose 
                                            $\varepsilon$ and $\delta$ small enough so that
                                            \begin{equation}
  \|\nabla\frak{q}\|_{L^2(\C)}\leq \sigma.
\end{equation}

The next  lemma shows that if $\sigma$ is small enough then \eqref{eq:gauge:2est}, \eqref{eq:gauge:pest}
                                            and \eqref{eq:gauge:ppest} 
 hold for a uniform constant $\Theta.$
\begin{Lm}\label{la:qgaugeuniformest}
There exists a $\Theta > 0$ and a $\sigma > 0$ so that whenever $\frak{q} \in \dot{W}^{1,p} \cap \dot{W}^{1,p^{\prime}}(\C,({\mathcal{U}}({\K}))$ and $\frak{q}-1 \in L^{\frac{2p}{p-2}}(\C,{\K})$ so that \eqref{decomp} is satisfied and it holds
\begin{equation}
  \|\nabla\frak{q}\|_{L^2(\C)}\leq \sigma,
\end{equation}
then  \eqref{eq:gauge:2est}, \eqref{eq:gauge:pest}
                                            and \eqref{eq:gauge:ppest} 
 hold.\hfill $\Box$
\end{Lm}
\noindent{\bf Proof of Lemma \ref{la:qgaugeuniformest}.}\par
Let us write   $\omega =\Delta\mu$ with $\nabla\mu\in L^{p}\cap L^{p^{\prime}}$. Let  $\xi\in \dot W^{1,p}\cap \dot W^{1,p^{\prime}}(\C,\R)$ be  such that
 \begin{equation}\label{estgauge7}
 \left\{\begin{array}{c}
 \Pi_i(\frak{q}^{-1}\partial_{x_1} \frak{q}- \p_{x_1}\mu )=-\partial_{x_2}\xi i\\[5mm]
 \Pi_i(\frak{q}^{-1}\partial_{x_2} \frak{q}- \p_{x_2}\mu )=\partial_{x_1}\xi i.\end{array}\right.
 \end{equation}
 Then
 \begin{equation}\label{xiequ}
 -\Delta \xi i= \Pi_i(\partial_{x_2}(\frak{q}^{-1}\partial_{x_1} \frak{q}))-\Pi_i(\partial_{x_1}(\frak{q}^{-1}\partial_{x_2} \frak{q})).
 \end{equation}
For every $q\in [p^{\prime},p]$ it holds \footnote{We use the fact that if $\nabla a\in L^{2,\infty}$, $\nabla b\in L^{q}$, with $q\in [p^{\prime},p]$ and  if
  $-\Delta \phi=\nabla a\cdot\nablap b,$ in  $\C$ then $\nabla \phi\in L^q$ with $\|\nabla\phi\|_{L^q}\leq C_q\|\nabla b\|_{L^q}\|\nabla a\|_{L^{2,\infty}}$. The constant $C_q$ is uniformly bounded if $q\in [p^{\prime},p]$ (see \cite{Hel}).
  }

\begin{eqnarray}\label{estgauge15}
 \|\nabla\xi\|_{L^q}&\lesssim &\|\nabla \frak{q}\|_{L^{2,\infty}}\|\nabla \frak{q}\|_{L^{q}}
 \lesssim \sigma \|\nabla \frak{q}\|_{L^q}
 \end{eqnarray}
We can write
\begin{eqnarray}
\frak{q}^{-1}\partial_{x_1} \frak{q}-\frak{q}^{-1}\partial_{x_2} \frak{q}i&=&\Pi_i(\frak{q}^{-1}\partial_{x_1} \frak{q})-\Pi_i(\frak{q}^{-1}\partial_{x_2} \frak{q})i\nonumber\\
&+&\Pi_{jk}(\frak{q}^{-1}\partial_{x_1} \frak{q})-\Pi_{jk}(\frak{q}^{-1}\partial_{x_2} \frak{q})i\nonumber\\&=&
(\partial_{x_1}\xi-\partial_{x_2}\xi i)+\Pi_i(\p_{x_1}\mu-\p_{x_2}\mu)+\frak{g}
\end{eqnarray}
or equivalently
\begin{equation}\label{estgauge9}
 \partial_{x_1} \frak{q}- \partial_{x_2} \frak{q}i=\frak{q}((\partial_{x_1}\xi-\partial_{x_2}\xi i)+\Pi_i(\p_{x_1}\mu-\p_{x_2}\mu)+\frak{g}) .
\end{equation}
 Therefore by combining \eqref{estgauge15} and \eqref {estgauge9}
 we get   for every $q\in [p^{\prime},p]$
  
     \begin{eqnarray}\label{estgauge16}
 \|\nabla \frak{q}\|_{L^q}&\le & C[ |\nabla \xi\|_{L^q}+ \|\nabla\mu\|_{L^q}+\|\frak{g}\|_{L^q}]\nonumber\\[5mm]
 &\le & C\sigma \|\nabla \frak{q}\|_{L^q}+C\|\omega\|_{ \dot{W}^{-1,q}}+C \|\frak{g}\|_{L^q}.
\end{eqnarray}
Choosing $\Theta=\frac{C}{1-C\sigma}$ we have 
$$
 \|\nabla \frak{q}\|_{L^q}\le \Theta( \|\omega\|_{ \dot{W}^{-1,q}}+\|\frak{g}\|_{L^q}).$$
 This concludes the proof of lemma~\ref{la:qgaugeuniformest}.\hfill $\Box$
 
 \medskip
 
 \noindent{\bf End of the proof of Proposition~\ref{pr:VepseqUeps}.}
 As we have already observed if $\varepsilon$ is small enough the fact that $\omega,\frak{g}\in {\mathcal{U}}_{\varepsilon}$ implies 
  that  $\frak{q}\in {\mathcal{V}}_{\varepsilon, \Theta}$, it satisfies 
 $ \|\nabla \frak{q}\|_{L^2}\le  \sigma $ where $\sigma$ is the constant appearing in Lemma \ref{la:qgaugeuniformest}. Therefore  thanks to lemma~\ref{la:qgaugeuniformest} the openness property (iv) is proven and Proposition~\ref{pr:VepseqUeps} is then established.}
\hfill $\Box$

\medskip

\noindent{\bf Proof of Theorem \ref{mainth}.}

Let  $\frak{f} $ solve  \eqref{est8} with $\frak{q}\in { N}^{-1}(0,-\partial_{x_1} \alpha j+\partial _{x_2} \alpha k)$ and
$\|\nabla \frak{q}\|_{L^2}\le \Theta\|\nabla\alpha\|_{L^2}.$   By  definition $\frak{q}$ satisfies  
\begin{equation} \label{est10}
\left\{\begin{array}{c}
  \Pi_i(\partial_{x_1}(\frak{q}^{-1}\partial_{x_1} \frak{q})+\partial_{x_2}(\frak{q}^{-1}\partial_{x_2} \frak{q}))=0 \\[5mm]
 \Pi_{jk}(\frak{q}^{-1}\partial_{x_1} \frak{q}-\frak{q}^{-1}\partial_{x_2} \frak{q}\ i)=-\partial_{x_1}\alpha j+\partial_{x_2}\alpha k\,.
 \end{array}
 \right.
   \end{equation}
We analyze the first equation in  \eqref{est10}.
\par
We have
\begin{equation}\left\{\begin{array}{c}
\Pi_i(\partial_{x_1}(\frak{q}^{-1}\partial_{x_1} \frak{q})+\partial_{x_2}(\frak{q}^{-1}\partial_{x_2} \frak{q}))=0 \\[5mm]
\Updownarrow\\[5mm]
 \partial_{x_1}\left(\Pi_i(\frak{q}^{-1}\partial_{x_1} \frak{q})\right)+\partial_{x_2}\left(\Pi_i(\frak{q}^{-1}\partial_{x_2} \frak{q})\right)=0\,.
 \end{array}\right.
 \end{equation}
 Therefore there exists $\zeta\in \dot W^{1,2}(\C,i\R)$ such that
 \begin{equation}\label{est12}
 \left\{\begin{array}{c}
 \Pi_i(\frak{q}^{-1}\partial_{x_1} \frak{q})=-\partial_{x_2}\zeta\\[5mm]
 \Pi_i(\frak{q}^{-1}\partial_{x_2} \frak{q})=\partial_{x_1}\zeta\,.
 \end{array}\right.
 \end{equation}
From \eqref{est12} it follows in particular that
\begin{eqnarray}\label{est13-aa}
-\Delta\zeta&=& \partial_{x_2}\left(\Pi_i(\frak{q}^{-1}\partial_{x_1} \frak{q})\right)- \partial_{x_1}\left(\Pi_i(\frak{q}^{-1}\partial_{x_2} \frak{q})\right)\nonumber\\[5mm]
&=&\Pi_i\left( \partial_{x_2}(\frak{q}^{-1}\partial_{x_1} \frak{q})- \partial_{x_1}(\frak{q}^{-1}\partial_{x_2} \frak{q})\right)\,.
\end{eqnarray}
The right hand side of \eqref{est13} is a sum of jacobians, hence it is in the Hardy space ${\mathcal{H}}^1(\C)$.
It follows in particular that $\nabla\zeta\in L^{2,1}(\C),$ with $$\|\nabla\zeta\|_{L^{2,1}}\lesssim \|\nabla \frak{q}\|^2_{L^2}\quad .$$
 We have
 \begin{eqnarray}\label{est13}
 \frak{q}^{-1}\partial_{x_1} \frak{q}-\frak{q}^{-1}\partial_{x_2} \frak{q}\ i&=&\Pi_i(\frak{q}^{-1}\partial_{x_1} \frak{q})+\Pi_{jk}(\frak{q}^{-1}\partial_{x_1} \frak{q})
 -\Pi_i(\frak{q}^{-1}\partial_{x_2} \frak{q})\ i-\Pi_{jk}(\frak{q}^{-1}\partial_{x_2} \frak{q})\ i\nonumber\\[5mm]
 &=&-\partial_{x_2}\zeta-\partial_{x_1}\zeta \ i+\Pi_{jk}(\frak{q}^{-1}\partial_{x_1} \frak{q} -\frak{q}^{-1}\partial_{x_2} \frak{q}\ i) \\[5mm]
 &=&-2i (\partial_z\zeta)-\partial_{x_1}\alpha\ j+\partial_{x_2}\alpha\, k\,.\nonumber
 \end{eqnarray}
In \eqref{est13} {  we use the fact that  $\underline{ {\Pi_{jk}(ai)=\Pi_{jk}(a)i}}$ for ${ {a\in\K.}}$}
By combining \eqref{est8}, \eqref{est10} and \eqref{est13} we get

 \begin{equation}\label{est14-ss}
  \partial_{x_1}[\frak{q} {\frak{f}}]-\partial_{x_2}[\frak{q} i\ {\frak{f}}]= -\,2\,\frak{q}\,i\,[\partial_{z} \zeta] \,{\frak{f}}\ .
\end{equation}
We set
\begin{equation}\label{est15-aa}
-\Delta A=2\,\frak{q}\,i\,[\partial_{z} \zeta] \frak{f}\,.
\end{equation}
Observe that 
\begin{eqnarray}\label{est16-aa}
\|\nabla A\|_{L^{2,\infty}}&\lesssim&\|\frak{q}\frak{f}\ \nabla\zeta\|_{L^1}\lesssim   \|\nabla\zeta\|_{L^{2,1}}\|\frak{q\ f}\|_{L^{2,\infty}}\\[5mm]
&\lesssim& \|\nabla \frak{q}\|^2_{L^{2}}\ \| \frak{f}\|_{L^{2,\infty}} \lesssim\varepsilon_0^2\ \| \frak{q\,f}\|_{L^{2,\infty}}.\nonumber
\end{eqnarray}
Since
$$\partial_{x_1}(\frak{q} {\frak{f}}-\partial_{x_1}A)-\partial_{x_2}(\frak{q}\ i \frak{f}+\partial_{x_2}A)=0\,,$$
there exists $B\in \dot W^{1,(2,\infty)}$ such that
\begin{equation}
\left\{\begin{array}{c}
\frak{q} {\frak{f}}-\partial_{x_1}A=-\partial_{x_2} B\\[5mm]
-\frak{q}i\ \frak{f}-\partial_{x_2}A=\partial_{x_1} B\,.\end{array}\right.
\end{equation}
Therefore
\begin{equation}\label{est15}
\left\{\begin{array}{c}
\frak{f}=\frak{q}^{-1}(\partial_{x_1}A-\partial_{x_2} B)\\[5mm]
\frak{f}=i\frak{q}^{-1}(\partial_{x_2}A+\partial_{x_1} B)\,.\end{array}\right.
\end{equation}
From \eqref{est15} it follows
\begin{equation}\label{est16-ss}
\left\{\begin{array}{l}
\partial_{x_1} A-\partial_{x_2} B=\frak{q}i\frak{q}^{-1}(\partial_{x_2}A+\partial_{x_1} B) \\[5mm]
-\partial_{x_2} A-\partial_{x_1} B=\frak{q}i\frak{q}^{-1}(\partial_{x_1}A-\partial_{x_2} B) \\[5mm]
-\partial_{x_2} B=\frak{q}i\frak{q}^{-1}(\partial_{x_2}A+\partial_{x_1} B)-\partial_{x_1} A \\[5mm]
-\partial_{x_1} B=\frak{q}i\frak{q}^{-1}(\partial_{x_1}A-\partial_{x_2} B)+\partial_{x_2} A\\[5mm]
-\Delta B=\partial_{x_1}\left(\frak{q}i\frak{q}^{-1}(\partial_{x_1}A-\partial_{x_2} B)\right)+\partial_{x_2}\left(\frak{q}i\frak{q}^{-1}(\partial_{x_2}A+\partial_{x_1} B)\right)\,.
\end{array}\right.
\end{equation}
We observe that $-\partial_{x_1}[\frak{q}i\frak{q}^{-1}\partial_{x_2} B]+\partial_{x_2}[\frak{q}i\frak{q}^{-1}\partial_{x_1} B]$ is sum of Jacobians and therefore we can apply 
Wente's Lemma (case $L^2-L^{2,\infty}$) :
\begin{eqnarray}\label{est16}
\|\nabla B\|_{L^{2,\infty}}&\lesssim& \|\nabla A\|_{L^{2,\infty}}+\|\nabla \frak{q}\|_{L^{2}}\|\nabla B\|_{L^{2,\infty}}\nonumber\\[5mm]
&\lesssim& \|\nabla A\|_{L^{2,\infty}}+\varepsilon_0\|\nabla B\|_{L^{2,\infty}}.
\end{eqnarray}
Estimate \eqref{est16} implies that
\begin{eqnarray}\label{est17-bb}
\|\nabla B\|_{L^{2,\infty}}&\lesssim& \|\nabla A\|_{L^{2,\infty}}\nonumber\\[5mm]
\|\frak{q} {\frak{f}}\|_{L^{2,\infty}} &\lesssim& \|\nabla A\|_{L^{2,\infty}}+ \|\nabla B\|_{L^{2,\infty}}\nonumber\\[5mm]
&\lesssim& \varepsilon^2_0\|\frak{q\ f}\|_{L^{2,\infty}}.
\end{eqnarray}
If $\varepsilon_0$ is small enough then $\frak{f}\equiv 0$. This concludes the proof of Theorem \ref{mainth}.~\hfill $\Box$
 \subsection{Morrey-Type Estimates}
 In this section we prove Morrey-type estimates for solutions to \eqref{maineq} in the case $n=2.$
 \begin{Prop}\label{regul}
 Let    $S\in W^{1,2}(\C,O(2))$   with $S^2=Id$ and $u\in L^{2}(\C)$ be a solution of \eqref{maineq}. Then $u\in W^{1,p}_{loc}$ for every $p\in [1,2).$ \hfill $\Box$
 \end{Prop}
\noindent{ \bf Proof of Proposition \ref{regul}.}
{\bf Step 1.} 
{  Assume that   $\|\nabla S\|_{L^2(B(0,1))}\leq \varepsilon_0$. 
\par
 {\bf Claim:}  There is $0<\varepsilon_0<1$ and  $\tilde S\in \dot W^{1,2}(\C,Sym(2))$ with $\tilde S^2=Id$ such that
$\tilde S=S$ in $B(0,1)$ and $\|\nabla \tilde S\|_{L^2(\R^2)}\le C \|\nabla S\|_{L^2(B(0,1))}$.
\par
 { For the proof of the claim we refer to \cite{Riviere}}.\par}
 
Now let $v\in L^{2}(\R^2)$ be such that
$$\nabla^{\perp} v=S\nabla u~~\mbox{in ${\mathcal{S}}^{\prime}(\R^2).$}$$
 
%Consider now $\chi\in C^{\infty}_c(\C)$ with $spt(\chi)\subseteq B(0,1).$ Then
%\begin{eqnarray}
%\nabla^{\perp}(\chi v)&=&\nabla^{\perp} \chi v+\chi\nabla^{\perp} v\\ &=&
 %\nabla^{\perp}\chi v +\chi S \nabla u=\underbrace{\nabla^{\perp}\chi v+\nabla(\chi S)u}_{\in L^1}-\underbrace{\nabla(\chi Su)}_{\in  \dot{W}^{-1,2}}\nonumber
%\end{eqnarray}
%Hence by applying Lemma \ref{BBL} we get $v\in L^2_{loc}(B(0,1)).$   

   By arguing as in the previous section we can find
 $\frak{q}\in \dot{W}^{1,2}(\C,{\mathcal{U}}(\K))$ with 
 $\|\nabla \frak{q}\|_{L^2(\C)}\le C \|\nabla \tilde S\|_{L^2(\C)}$ and $\zeta\in \dot{W}^{1,(2,1)}(\C)$ with $\|\partial_{z} \zeta\|_{L^{2,1}(\C)}\le \varepsilon^2_0$ such that
  \begin{equation}\label{est14}
   \partial_{x_1}[\frak{q} {\frak{f}}]-\partial_{x_2}[\frak{q} i\ {\frak{f}}]= -2\frak{q}\,i\,[\partial_{z} \zeta] {\frak{f}} ~~\mbox{in ${\mathcal{D}}^{\prime}(B(0,1))$} 
\end{equation}
 1. Let $x\in B(0,1/2)$ and $0<r<1/4$. We consider
\begin{equation}\label{estA}
\left\{\begin{array}{cc}
-\Delta A=2\,\frak{q}\,i\,[\partial_{z} \zeta]\ \frak{f}\,~~& \mbox{in $B(x,r)$}\\
A=0 ~~& \mbox{on $\partial B(x,r)$}.\end{array}\right.
\end{equation}
 We have
 \begin{equation}\label{estA2}
 \|\nabla A\|_{L^{2,\infty}(B(x,r))}\lesssim \varepsilon_0^2\ \|f\|_{L^{2,\infty}(B(x,r))}\,.
 \end{equation}
2.    Since
$$\partial_{x_1}(\frak{q} {\frak{f}}-\partial_{x_1}A)-\partial_{x_2}(\frak{q}\ i \frak{f}+\partial_{x_2}A)=0\,,$$
there exists $B\in \dot W^{1,(2,\infty)}(B(x,r))$ such that
\begin{equation}\label{ABf}
\left\{\begin{array}{c}
\frak{q} {\frak{f}}-\partial_{x_1}A=-\partial_{x_2} B\\[5mm]
-\frak{q}i\ \frak{f}-\partial_{x_2}A=\partial_{x_1} B\,.\end{array}\right.
\end{equation}
 We have 
\begin{equation}\label{DeltaB}
-\Delta B=\partial_{x_1}\left(\frak{q}i\frak{q}^{-1}(\partial_{x_1}A-\partial_{x_2} B)\right)+\partial_{x_2}\left(\frak{q}i\frak{q}^{-1}(\partial_{x_2}A+\partial_{x_1} B)\right)~~\mbox{in ${\mathcal{D}}^{\prime}(B(x,r))$}
\end{equation}
 
We decompose $B=\beta_1+\beta_2$ in $B(x,r)$ where $\beta_1$ and $\beta_2$ satisfy respectively
 \begin{equation}\label{systB}
 \left\{\begin{array}{cc}
\ds \Delta \beta_1=0&~~\mbox{in $B(0,r)$}\\[5mm]
\ds \beta_1=B& ~~\mbox{on $\partial B(0,r)$}\end{array}\right. 
  ~~\mbox{and}~~~ \left\{\begin{array}{cc}
\ds \Delta \beta_2=\Delta B&~~\mbox{in $B(0,r)$}\\[5mm]
 \beta_2=0 &~~\mbox{on $\partial B(0,r)$}\end{array}\right.\end{equation}
The following estimates hold:

\medskip

\noindent{\bf Estimate of $\beta_2$:} \par
Wente inequality ($L^{2,\infty}-L^2$ case) combined with classical Calderon Zygmund inequalities give
\begin{equation}\label{estimatebeta2}
\|\nabla\beta_2\|_{L^{2,\infty}(B(x,r))}\le \|\nabla A\|_{L^{2,\infty}(B(x,r))}+\varepsilon_0\|\nabla B\|_{L^{2,\infty}(B(x,r))}\,.
\end{equation}
\noindent{\bf Estimate of $\beta_1$:} \par
Since $\beta_1$ is harmonic, for every $0<\delta<\frac{3}{4}$  we have
 \begin{eqnarray}\label{estbeta1}
 \|\nabla\beta_1\|^2_{L^{2,\infty}(B(x,\delta r))}&\le& \|\nabla\beta_1\|^2_{L^2(B(x,\delta r))} \\
 &\le& \left(\frac{4\delta}{3}\right)^2  \|\nabla\beta_1\|^2_{L^2(B(x,3/4 r))}\le C\left(\frac{4\delta}{3}\right)^2  \|\nabla\beta_1\|^2_{L^{2,\infty}(B(x, r))},\nonumber
 \end{eqnarray}
 where $C$ is a constant independent of $r.$ In \eqref{estbeta1} we use the fact that the $L^{2,\infty}$ of the gradient of a harmonic function on the ball
 $B(x, r)$ controls all its other norms in balls $B(x,\eta r)$ with $\eta< 3/4.$

 \medskip
 
\noindent{\bf Estimate of $B$:}\par

Combining the previous estimates we obtain
 \begin{eqnarray}\label{estB}
 \|\nabla B\|_{L^{2,\infty}(B(x,\delta r))}&\lesssim & \|\nabla\beta_1\|_{L^{2,\infty}(B(x,\delta r)}+\|\nabla\beta_2\|_{L^{2,\infty}(B(x,\delta r)}\nonumber\\
 &\lesssim & \left(\frac{4\delta}{3}\right)  \|\nabla\beta_1\|_{L^{2,\infty}(B(x, r)}+  \|\nabla A\|_{L^{2,\infty}(B(x,r))}+\varepsilon_0\|\nabla B\|_{L^{2,\infty}(B(x,r))}\nonumber\\
   &\lesssim& \left(\frac{4\delta}{3}\right)\left[\|\nabla\beta_2\|_{L^{2,\infty}(B(x, r)}+\|\nabla B\|_{L^{2,\infty}(B(x,r))}\right]\nonumber\\
&+&   \|\nabla A\|_{L^{2,\infty}(B(x, r))}+\varepsilon_0\|\nabla B\|_{L^{2,\infty}(B(x,r))}
    \nonumber\\
 &\lesssim& 
  \left(\frac{4\delta}{3}\right)\left[ \|\nabla A\|_{L^{2,\infty}(B(x,r))}+\varepsilon_0\|\nabla B\|_{L^{2,\infty}(B(x,r))}+\|\nabla B\|_{L^{2,\infty}(B(x,r))}\right]\nonumber\\
  &+&   \|\nabla A\|_{L^{2,\infty}(B(x,  r))}+\varepsilon_0\|\nabla B\|_{L^{2,\infty}(B(x,  r))}
\nonumber\\
&\lesssim&\left[\left(\frac{4\delta}{3}\right)\varepsilon_0^2+\varepsilon_0^2\right] \|\frak{f}\|_{L^{2,\infty}(B(0,r))} \\&+&
\left[\left(\frac{4\delta}{3}\right)\varepsilon_0+\left(\frac{4\delta}{3}\right)+\varepsilon_0\right] \|\nabla B\|_{L^{2,\infty}(B(x,r))}\nonumber
  \end{eqnarray}
 Since $\|\nabla B\|_{L^{2,\infty}(B(x,r))}\le  \|\nabla A\|_{L^{2,\infty}(B(x,r))}+ \|\frak{f}\|_{L^{2,\infty}(B(x,r))}$ from \eqref{estB} one deduces that
  \begin{eqnarray}\label{estBbis}
 \|\nabla B\|_{L^{2,\infty}(B(x,\delta r))}&\lesssim & \left[\left(\frac{4\delta}{3}\right)\varepsilon_0^2+\varepsilon_0^2\right] \|\frak{f}\|_{L^{2,\infty}(B(0,r))}\\&+&
 \left[\left(\frac{4\delta}{3}\right)\varepsilon_0+\left(\frac{4\delta}{3}+\varepsilon_0\right)\right](1+\varepsilon_0^2)\|\frak{f}\|_{L^{2,\infty}(B(0,r))}.\nonumber
 \end{eqnarray}
By combining  \eqref{ABf} and \eqref{estBbis} we obtain
  \begin{eqnarray}\label{estf}
  \|\frak{f}\|_{L^{2,\infty}(B(x,\delta r))}&\lesssim & \|\nabla A\|_{L^{2,\infty}(B(x,\delta r))}+ \|\nabla B\|_{L^{2,\infty}(B(x,\delta r))}\nonumber\\
  &\lesssim& \gamma \|\frak{f}\|_{L^{2,\infty}(B(0,r))}\,.
  \end{eqnarray}
  where $\gamma=\gamma(\delta,\varepsilon_0)<1.$
  By iterating \eqref{estf} we get the existence of a constant $0<\alpha<1$ such that
  \begin{equation}\label{estf2}
  \sup_{x\in B(0,1/2),0<r<1/4}r^{-\alpha}\|f\|_{L^{2,\infty}(B(x,  r))}<+\infty\,.
  \end{equation}
 Now we plug the estimate \eqref{estf2} into \eqref{estA} and we get
    \begin{equation}\label{estA3}
  \sup_{x\in B(0,1/2),0<r<1/4}r^{-\alpha}\|\Delta A\|_{L^{1}(B(x, r))}<+\infty\,.
  \end{equation}
  and  therefore
  \begin{equation}\label{estA4}
  \sup_{x\in B(0,1/2),0<r<1/4}r^{-\alpha}\|\nabla A\|_{L^{2,\infty}(B(x, r))}<+\infty\,.
  \end{equation}
From \eqref{estA3} it follows in particular that $\nabla A\in L^{q}(B(0,1/4))$ for all $q<\frac{2-\alpha}{1-\alpha}$ (See again Adams \cite{Ad}, Remark after Proposition 3.2).
\par
From \eqref{ABf}, \eqref{estf2}, \eqref{estA4} it follows that
\begin{equation}\label{estB2}
  \sup_{x\in B(0,1/2),0<r<1/4}r^{-\alpha}\|\nabla B\|_{L^{2,\infty}(B(x,  r))}<+\infty.
  \end{equation}
 By plugging \eqref{estB2} into \eqref{DeltaB} and \eqref{systB} one gets that
 \begin{equation}\label{estB3}
  \sup_{x\in B(0,1/2),0<r<1/4}r^{-\alpha}\|\Delta B\|_{L^{2,\infty}(B(x,  r))}<+\infty.
  \end{equation}
which implies that $\nabla B\in L^{q}(B(0,1/4))$ for all $q<\frac{2-\alpha}{1-\alpha}$ as well. Therefore $\frak{f}\in  L^{q}(B(0,1/4))$ for all $q<\frac{2-\alpha}{1-\alpha}$ as well. Actually one can show by bootstrap arguments that  $\frak{f}\in  L^{q}_{loc}$ for all $q<+\infty.$
 
 \medskip
 
 \noindent{\bf Step 2.} From Step 1 it follows that $Su\in L^{q}_{loc}(\C)$ for all $q<+\infty$.
 Since $u$ solves \eqref{maineq} we have
 \begin{equation}\label{maineqbis}
 \Delta (Su)=\div(\nabla(Su))=\div(\nabla SS \ Su) ~~\mbox{in ${\mathcal{D}}^{\prime}(\C).$}
 \end{equation}
  From \eqref{maineqbis} one gets that $\nabla(Su)\in L^{\frac{2q}{q+2}}_{loc}$ for all $q<+\infty$ and therefore
  $\nabla u=\nabla S (Su)+S\nabla (Su) \in L^{p}_{loc}$ for all $p<2.$
  This concludes the proof of proposition~\ref{regul} which itself implies theorem~\ref{th-intro-1} in the case of $2$-D codomains.   \hfill$\Box$

  \section{Proof of theorem~\ref{th-intro-1} : the   general case $n\ge 2$}
We are going to present here another approach  to study the regularity of the equation \eqref{maineq} which works for every $n\ge 2$. We start by showing the bootstrap test:
\begin{Th}\label{bootstrap2}
Let  $S\in \dot W^{1,2}(\C, O(n))$ with $S^2=Id $   and $u\in L^2(\C,\R^n)$ be a solution to the equation \eqref{maineq}. There is $\varepsilon_0>0$ such that
if   $\|\nabla S\|_{L^2(\C)}\le \varepsilon_0, $
then $u\equiv 0.$  
\end{Th}
\noindent{\bf Proof of theorem \ref{bootstrap2}}
From Lemma \ref{BBL} we can find  
 $v\in L^{2}(\C,\R^n)$   such that $\nablap v= S\,\nabla u.$\par
Assume that $\|\nabla S\|_{L^2(\C)}\le \varepsilon_0 $ where $\varepsilon_0$ is the constant appearing in Theorem \ref{S}. Then  there is    $Q\in \dot{W}^{1,2}(\C,SO(n))$ such that
$$S=Q^{-1}\,S^0\,Q$$
where  $S^0$ is the matrix \eqref{S0g} and $\|\nabla Q\|_{L^2(\C)}\lesssim \varepsilon_0.$\par

 We set 
 $$
f=f_{\Re}+if_{\Im} =S^0Qu+i Qv.$$
Equation \eqref{maineq-ff} is equivalent to the system:
\begin{equation}
\left\{\begin{array}{c}
\partial_{x_1} f_{\Re}+\partial_{x_2}f_{\Im}=S^0\partial_{x_1}QQ^{-1}S^0f_{\Re}+\partial_{x_2}QQ^{-1}f_{\Im}\\[5mm]
-\partial_{x_2} f_{\Re}+\partial_{x_1}f_{\Im}=-S^0\partial_{x_2}QQ^{-1}S^0f_{\Re}+\partial_{x_1}QQ^{-1}f_{\Im}\,.
\end{array}
\right.
\end{equation}
 We can write 
 $$S^0=((-1)^{\min(2m+1,2i)} \delta_{ij})_{1\le i,j\le n}\,.$$
Let $\Omega=(\omega_{ij})_{1\le i, j\le n}$ be an anti-symmetric real matrix (i.e. $\omega_{ij}=-\omega_{ji}$), then
$$\tilde\Omega=S^0\,\Omega\, S^0=\left(\omega_{ij}(-1)^{\min(2m+1,2i)+\min(2m+1,2j )}\right)\,.$$
Therefore 
\begin{equation}\label{link}
\left\{\begin{array}{c}
  \tilde{\omega}_{ij}=\omega_{ij} ~~\Longleftrightarrow ~~\mbox{$i,j\le m$ and $i,j>m$}\\[5mm]
 \tilde{\omega}_{ij}=-\omega_{ij} ~~\Longleftrightarrow ~~\mbox{otherwise}\,.
\end{array}
\right.
\end{equation}
Observe that the matrix $\ti{\Om}$ is still anti-symmetric. We set $\Omega^\ell:=\partial_{x_\ell}QQ^{-1}$ and $\tilde\Omega^{\ell}= S^0\,\partial_{x_\ell}QQ^{-1}\,S^0.$
\begin{equation}
\left\{\begin{array}{c}
\partial_{x_1} f_{\Re}+\partial_{x_2}f_{\Im}= \tilde\Omega^1\ f_{\Re}+\Omega^2\ f_{\Im} \\[5mm]
\partial_{x_2} f_{\Re}-\partial_{x_1}f_{\Im}=\tilde\Omega^2\ f_{\Re}-\Omega^1\ f_{\Im}\,.
\end{array}
\right.
\end{equation}
Then we get
\begin{eqnarray}\label{estnb2-aa}
(\partial_{x_1}-i\partial_{x_2})(f_{\Re}+if_{\Im})&=& \tilde\Omega^1\ f_{\Re}+\Omega^2\ f_{\Im}-i(\tilde\Omega^2\ f_{\Re}-\Omega^1\ f_{\Im})\nonumber\\[5mm]
&=&\left(\tilde\Omega^1-i\tilde\Omega^2\right)\ f_{\Re}+i\left( \Omega^1-i\Omega^2\right)\ f_{\Im}\nonumber\\[5mm]
&=&\left(\frac{\tilde\Omega^1-i\tilde\Omega^2}{2}\right)\ ((f_{\Re}+i f_{\Im})+(f_{\Re}-i f_{\Im}))\nonumber\\[5mm]
 &+&\left(\frac{\Omega^1-i\Omega^2}{2}\right)\ ((f_{\Re}+i f_{\Im})-(f_{\Re}-i f_{\Im}))\,.
 \end{eqnarray}
 Which gives
 \begin{eqnarray}\label{estnb2}
(\partial_{x_1}-i\partial_{x_2})(f_{\Re}+if_{\Im})&=& \frac{1}{2}\left[(\tilde\Omega^1+\Omega^1)-i(\tilde\Omega^2+\Omega^2)\right](f_{\Re}+i f_{\Im})\nonumber\\[5mm]
&&+ \frac{1}{2}\left[(\tilde\Omega^1-\Omega^1)-i(\tilde\Omega^2-\Omega^2)\right](f_{\Re}-i f_{\Im})\,.
 \end{eqnarray}
 From \eqref{link} it follows for $\ell=1,2$
 \begin{equation}
\frac{\tilde\Omega^\ell+\Omega^\ell}{2}=\left(\begin{array}{c|c}
\omega_{ij}^{\ell}&0_{m\times n-m}\\[5mm] \hline &\\
0_{n-m\times m} & \omega_{ij}^{\ell}  \end{array}\right) 
\end{equation}
 and
 \begin{equation}
\frac{\tilde\Omega^\ell-\Omega^\ell}{2}=\left(\begin{array}{c|c}
0_{m\times m} &-\omega_{ij}^{\ell} \\[5mm] \hline &\\
-\omega_{ij}^{\ell} & 0_{n-m\times n- m}  \end{array}\right) \,.
\end{equation}
We can write the system \eqref{estnb2} as
\begin{equation}\label{estnb21}
\partial_L f=\frac{1}{2}\Omega^{+}f+\frac{1}{2}\Omega^-\bar f\,,
\end{equation}
where
\begin{eqnarray} 
\Omega^{+}&=&\frac{(\tilde\Omega^1+\Omega^1)-i(\tilde\Omega^2+\Omega^2)}{2}\label{Omega1}\\
\Omega^{-}&=&\frac{(\tilde\Omega^1-\Omega^1)-i(\tilde\Omega^2-\Omega^2)}{2}\label{Omega2}\,.
\end{eqnarray}
We observe that by construction for every $i,j$ we have
\begin{eqnarray} 
\Im(\partial_{\bar z}\Omega^{+}_{ij})&=&\partial_{x_2}(\Omega^{+}_{ij} )^{\Re}+\partial_{x_1}(\Omega^{+}_{ij})^{\Im}\in {\mathcal{H}}^1(\R^2)\label{jacproperty1}\\ 
\Im(\partial_{\bar z}\Omega^{-}_{ij})&=&\partial_{x_2}(\Omega^{-}_{ij} )^{\Re}+\partial_{x_1}(\Omega^{-}_{ij})^{\Im}\in {\mathcal{H}}^1(\R^2)\label{jacproperty2}
\end{eqnarray}
with 
$$\|\Im(\partial_{\bar z}\Omega^{+}_{ij})\|_{{\mathcal{H}}^1(\R^2)}\lesssim\|\nabla Q\|^2_{L^2(\R^2)}, \|\Im(\partial_{\bar z}\Omega^{-}_{ij})\|\lesssim\|\nabla Q\|^2_{L^2(\R^2)},$$
since these quantities are linear combinations of Jacobians of functions  (the components of the matrix $Q$) with gradient in $L^2$.

\medskip

Let $M$ be defined as follows:
\begin{equation}\label{matrixQ}
M=\left(\begin{array}{c|c}
M_R  &0_{m\times n-m}\\\hline
0_{n-m\times m} & M_L\end{array}\right) \end{equation}
where $M_R\in \dot W^{1,2}(\R^2,SO(m))$ and $M_L\in  \dot W^{1,2}(\R^2,SO(n-m))$.
The following identity holds
\begin{eqnarray}\label{systemQf}
\partial_L(Mf)&=&\partial_L M f+M\partial_z f=(\partial_L MM^{-1})(Mf)+\frac{1}{2}M\left(\Omega^{+}f+\Omega^{-}\bar f\right)\nonumber\\
&=&\left(\partial_L M M^{-1}+\frac{1}{2}M\Omega^{+}M^{-1}\right)Mf+\frac{1}{2}(M \Omega^{-}M^{-1})\overline{M f}\,.
\end{eqnarray}
\noindent{\bf  Claim 1:}  There are two constants $\varepsilon(n)>0$ and $C(n)>0$ depending only on $n$ such that if     $\|\Omega^+\|_{L^2}< \varepsilon(n)$  there exists  a matrix $M$ of the form \eqref{matrixQ} and $\eta\in \dot W^{1,(2,1)}(\R^2)\cap L^{\infty}(\R^2)$ such that  $$\partial_L M M^{-1}+\frac{1}{2}M\Omega^{+}M^{-1}=-i\partial_z\eta\,$$
an
$$\|\nabla M\|_{L^2(\R^2)},\, \|\nabla \eta\|_{L^{2,1}(\R^2)}\lesssim \|\Omega^+\|_{L^2}\ .$$

{\bf Proof of Claim 1.} By the same arguments in Lemma A.3 in  \cite{Riv1}  we can find
$M\in \dot W^{1,2}({\C},SO(n))$ of the form \eqref{matrixQ} with 
$\|\nabla M\|_{L^2(\R^2)}\lesssim \|\Omega^+\|_{L^2}$  and $\eta\in \dot W^{1,2}({\C},so(n))$ such that
 \be
 \lf\{
 \begin{array}{l}
\ds -\partial_{x_2}\eta =\partial_{x_1} M M^{-1}+ M[\tilde\Omega_1+\Omega_1]M^{-1}\\[5mm]
\ds \partial_{x_1}\eta =\partial_{x_2} M M^{-1}+ M[\tilde\Omega_2+\Omega_2]M^{-1}\,.
\end{array}
\rg.
 \ee
  It follows that
 \begin{eqnarray}\label{esteta-bb}
 -\Delta \eta&=&\underbrace{\partial_{x_2}(\partial_{x_1} M M^{-1})-\partial_{x_1}(\partial_{x_2} M M^{-1})}_{(1)}\\[5mm]
 &+& \underbrace{\partial_{x_2}\left(M\left(\tilde\Omega_1+\Omega_1\right)M^{-1}\right)-\partial_{x_1}\left(M\left(\tilde\Omega_2+\Omega_2\right)M^{-1}\right)}_{(2)}.
 \end{eqnarray}
 The first term $(1)$ on the right hand side of \eqref{esteta-bb} is in the Hardy Space ${\mathcal{H}}^1(\R^2)$ since it is a linear combination of Jacobians of
 functions with derivative in $L^2$.
 
\par
 {\bf Claim 2:} 
The second term   $(2)$  is in $ \dot{W}^{-1,(2,1)}(\R^2).$\par
{\bf Proof of the Claim 2.}  Indeed we observe that each component of $(2)$ can be written in the form
 \begin{equation}\label{esteta2}
 \partial_{x_2}(a\omega_1)-\partial_{x_1}(a\omega_2)
 \end{equation}

 where $a\in (\dot W^{1,2}\cap L^{\infty})(\R^2) $ and $\omega_1,\omega_2\in L^{2}(\R^2)$ satisfy
 \be
 \label{esteta2-ba}
 \partial_{x_2}\omega_1-\partial_{x_1}\omega_2\in  {\mathcal{H}}^1(\R^2)\,.
 \ee
 
 Actually we have
 \begin{eqnarray}\label{MOM}
 &&\partial_{x_2}\left(M\left(\tilde\Omega_1+\Omega_1\right)M^{-1}\right)_{ij}-\partial_{x_1}\left(M\left(\tilde\Omega_2+\Omega_2\right)M^{-1}\right)_{ij} \nonumber\\
  &=&\partial_{x_2}(M_{ik}(\tilde\Omega_1+\Omega_1)_{kt}M^{-1}_{tj})- \partial_{x_1}(M_{ik}(\tilde\Omega_2+\Omega_2)_{kt}M^{-1}_{tj}).
  \end{eqnarray}
  One sets
  $$a=M_{ik}M^{-1}_{tj},~~\omega_1=(\tilde\Omega_1+\Omega_1)_{kt},~~\omega_2= (\tilde\Omega_2+\Omega_2)_{kt}.$$
 Let $c,b\in   \dot{ W}^{1,2}(\R^2)$ be such that
 \[
 \lf(\begin{array}{c}\omega_1\\[3mm]
 \omega_2\end{array}\rg)=\nabla^{\perp} c+\nabla b\,.
 \]
 We can deduce from \eqref{esteta2-ba} that $\Delta c\in {\mathcal{H}}^1(\R^2)$ hence $c\in  \dot W^{1,(2,1)}(\R^2)$. We can now rewrite \eqref{esteta2} as follows
 \begin{eqnarray}\label{esteta3} 
 \partial_{x_2}[a(\partial_{x_1}b-\partial_{x_2}c)]- \partial_{x_1}[a(\partial_{x_2}b+\partial_{x_1}c)]&=&
 \partial_{x_2}a\partial_{x_1} b- \partial_{x_1}a\partial_{x_2} b\nonumber\\[5mm]&-&(\partial_{x_1}[a\partial_{x_1}c]+\partial_{x_2}[a\partial_{x_2}c]).
 \end{eqnarray}
 We observe that 
 $\partial_{x_2}a\partial_{x_1} b- \partial_{x_1}a\partial_{x_2} b\in {\mathcal{H}}^1(\R^2)$ and
 $\partial_{x_1}[a\partial_{x_1}c]+\partial_{x_2}[a\partial_{x_2}c]\in  \dot{W}^{-1,(2,1)}(\R^2)$.
 This gives that $(2)$ is in $ \dot{W}^{-1,(2,1)}(\R^2)$ and {\bf this concludes the proof of Claim 1 and Claim 2.}~$\Box$.

 \bigskip
 
 The system \eqref{systemQf} can then be written as
 \begin{equation}\label{systemQf2-ss}
  \partial_L(Mf)=A\, ( M f)+ B\,\overline{ M f}
  \end{equation}
  with $A=-i\partial_z\eta\in L^{2,1}(\R^2,{\mathcal{M}}_{n\times n}(\C))$ and $B= \frac{1}{2}(M \Omega^{-}M^{-1})\in L^{2}(\R^2,{\mathcal{M}}_{n\times n}(\C))$ with
  $$\|\nabla A\|_{L^{2,1}(\R^2)},\,\|\nabla B\|_{L^2(\R^2)}\lesssim \|\nabla M\|_{L^2}^2.$$
   $B$  satisfies
  for every $i,j$ ${B}_{ij}=-B_{ji}$ and 
  \begin{equation}\label{Bcondition}
  \partial_{x_2}(B^{\Re}_{ij})+\partial_{x_1}(B^{\Im}_{ij})\in  \dot{W}^{-1,(2,1)}(\R^2).
  \end{equation}
  {\bf Proof of the Claim 3.}
 For every $ij$ we have $B_{ij}= \frac{1}{2}(M_{ik}\Omega^{-}_{kt}M^{-1}_{tj})$. We know from \eqref{jacproperty2}  that 
  \begin{eqnarray}  
\Im(\partial_{\bar z}\Omega^{-}_{kt})&=&\partial_{x_2}(\Omega^{-}_{kt} )^{\Re}+\partial_{x_1}(\Omega^{-}_{kt})^{\Im}\in {\mathcal{H}}^1(\R^2)\label{jacproperty3}
\end{eqnarray}
\begin{equation}\label{estomegaminus}
\|\Im(\partial_{\bar z}\Omega^{-}_{kt})\|_{{\mathcal{H}}^1(\R^2)}\lesssim\|\nabla Q\|^2_{L^2(\R^2)} \end{equation}
 
We proceed as in the proof of Claim 2:
let $c,b\in   \dot{ W}^{1,2}(\R^2)$ be such that
 \[
 \lf(\begin{array}{c}(\Omega^{-}_{kt} )^{\Re}\\[3mm]
 -(\Omega^{-}_{kt} )^{\Im}\end{array}\rg)=\nabla^{\perp} c+\nabla b\,.
 \]
 We can deduce from \eqref{estomegaminus} that $\Delta c\in {\mathcal{H}}^1(\R^2)$ hence $c\in  \dot W^{1,(2,1)}(\R^2)$. 
 Then setting $a=M_{ik}M^{-1}_{tj}$ we have
 \begin{eqnarray}
 &&\partial_{x_2}(a\Omega^{-}_{kt} )^{\Re})+\partial_{x_1}(a\Omega^{-}_{kt} )^{\Im}=
  \partial_{x_2}(a(\partial_{x_1}b-\partial_{x_2}c)+\partial_{x_1}(a(-\partial_{x_2}b-\partial_{x_1}c))\nonumber\\
  &=& \partial_{x_2} a \partial_{x_1}b-\partial_{x_1}a\partial_{x_2}b-(\partial_{x_2}(a\partial_{x_2}c)+\partial_{x_1}(a\partial_{x_1}c).
  \end{eqnarray}
  We then conclude as in the proof of Claim 3.~~$\Box$
  \par
  \medskip
  In the sequel we can focus our attention  to a  system of the type:
     \begin{equation}\label{system M f2}
  \partial_L\, g= A\,g+ B\,\overline{g}
  \end{equation}
  where $A\in L^{2,1}(\R^2,{\mathcal{M}}_{n\times n}(\C))$ and $B\in L^{2}(\R^2,{\mathcal{M}}_{n\times n}(\C))$  satisfying ${B}_{ij}=-B_{ji}$ and \eqref{Bcondition}.
  
  \medskip
  
 \noindent {\bf Step 1.} We first observe that
   \begin{equation}\label{systemQf3}
  \partial_L\, g= Ag -Bj\ g\,j
  \end{equation}
  where $j$ is the quaternion number satisfying $j^2=-1$ and $ij=-ji$.
  
  \medskip
  
 \noindent {\bf Step 2.} The function $gj$ satisfies the system
  \begin{equation}\label{systemQf4}
  \partial_L\, gj = Agj+ Bj\ g.
  \end{equation}
  {\bf Step 3.} We set 
 \be
  \label{def-G}
  G=\begin{pmatrix} g^1\\ \vdots\\g^n\\ g^1j\\ \vdots \\ g^nj\end{pmatrix}
  \ee
  $G$ satisfies
  \begin{equation}\label{systemQf5}
  \partial_L G= \Gamma\, G +\Gamma_1 \,G\,,
\end{equation}
where
\begin{equation}\label{defOmega1}
\Gamma_1=\Gamma_1^{\Re}-i\Gamma_1^{\Im}=\left(\begin{array}{c|c}
A & 0_{n\times n}\\[2mm] \hline &\\
 0_{n\times n}&A \end{array}\right)\,,
\end{equation}
and
\begin{equation}\label{defOmega}
\Gamma=\Gamma^{\Re}+i\Gamma^{\Im}=\left(\begin{array}{c|c}
0_{n\times n} &-Bj\\[2mm] \hline &\\
B j&0_{n\times n}  \end{array}\right)\,,
\end{equation}
where we have set $B=B^{\Re}+iB^{\Im}$ and 
$$
\Gamma^{\Re}:=\left(\begin{array}{c|c}
0_{n\times n} &(-B)^{\Re}j\\[5mm] \hline &\\
B^{\Re}j&0_{n\times n}  \end{array}\right) 
\quad\quad\mbox{and}\quad\quad
\Gamma^{\Im}:=\left(\begin{array}{c|c}
0_{n\times n} &(-B)^{\Im}j\\[5mm] \hline &\\
B^{\Im}j&0_{n\times n}  \end{array}\right) \,.
$$
Observe that
 $$
\ov{\Gamma^{\Re}}:=\left(\begin{array}{c|c}
0_{n\times n} &B^{\Re}j\\[5mm] \hline &\\
-B^{\Re}j&0_{n\times n}  \end{array}\right) 
$$
and then
 $$
(\ov{\Gamma^{\Re}})^t=\left(\begin{array}{c|c}
0_{n\times n} &-(B^t)^{\Re}j\\[5mm] \hline &\\
(B^t)^{\Re}j&0_{n\times n}  \end{array}\right) =-\Gamma^{\Re}
$$
Therefore
\begin{equation}\label{OmegaR}
(\ov{\Gamma^{\Re}})^t+\Gamma^{\Re}=0.
\end{equation}
Similarly  $(\ov{\Gamma^{\Im}})^t+\Gamma^{\Im}=0.$\par
  
Since   the coefficients of $\Gamma^{\Im}$ are in $j\,{\R}$, we obtain $\ov{i \Gamma^{\Im}}=-\ov{\Gamma^{\Im}}\,i=i\,\ov{\Gamma^{\Im}}$. Hence finally we have established $$(\ov{i \Gamma^{\Im}})^t=i\,(\ov{\Gamma^{\Im}} )^t=-\,i\,\Gamma^{\Im} $$\par\bigskip
The matrix $\Gamma=\Gamma^\Re+i\Gamma^\Im$  satisfies then $(\ov{\Gamma})^t+\Gamma=0$ which means that it belongs to the Lie algebra $u(2n,\K)$ of the {\it hyper-unitary group} $U(2n,\K)$. This is  the compact Lie group  of invertible $2n\times 2n$ quaternions matrices 
$D$ satisfying $\bar D^t D=D\bar D^t=Id_n.$ 
 
 \par
 \medskip
{ Let $G$ be a $L^2$ solution of \eqref{systemQf5} with $\Gamma\in L^2(\R^2,u(2n,\K)),\Gamma_1\in L^{2,1}(\R^2,{\mathcal{M}}_{2n\times 2n}(\K))$.  Let  us take $P\in L^2(\R^2,U(n,\K))$ (to be fixed later), then
the following estimates hold
 \begin{eqnarray}
&&\partial_{x_1}(PG)-\partial_{x_2}(PiG)=P\left[(\partial_{x_1}G-i\partial_{x_2} G)+P^{-1}\partial_{x_1}P -P^{-1}\partial_{x_2}P i \right]G\nonumber\\[5mm]
&=&P\left[2(\Gamma_1+\Gamma )+(P^{-1}\partial_{x_1}P -P^{-1}\partial_{x_2}P i)\right]P^{-1}(PG).
 \label{newsystem2}
\end{eqnarray}
 The key point is to choose $P$  in order to absorb in \eqref{newsystem2}    the term $2\Gamma$. \par
We first observe that if $P\in U(2n,\K)$ then $P^{-1}\nabla P\in u(2n,\K)$.  Actually since $P^{-1}=\overline{P}^t$ and $P^{-1}\nabla P=-\nabla{P^{-1}}P$ one has
$$\overline{(P^{-1}\nabla P)}^t=(\nabla \overline{P}^t)\overline{P^{-1}}^t=\nabla{P^{-1}}P=-P^{-1}\nabla P.$$ 
\footnote{  We recall that  the standard Hermitian form in $\K^n$ is defined  by $\langle x,  y\rangle:=\sum_{i=1}^n\bar x_i y_i$.
Therefore given $A,B$ two    $n\times n$ matrices with   entries in $\K$  we have
$$\langle AB x,y  \rangle=\langle B x,\bar A^ty  \rangle=\langle  x,\bar B^t\bar A^t y  \rangle$$
and therefore $\overline{(AB)}^t=\bar B^t\bar A^t,$ (see e.g. \cite{Z}, Section 4). }
We also recall that every matrix $U\in  u(2n,\K)$ can be represented as
$$
U=U_0+U_1i+U_2j+U_3k,$$
where $U_i$ are real $2n\times 2n$ matrices such that $U_0^t=-U_0$ and $U_i^t=U_i.$ for $i=1,2,3.$\par
Now we are going to proceed as in Section \ref{boot}.\par 
In the sequel we will denote by $\mbox{MSpan$\{1,i\}$}$ the space of $n\times n$ matrices  
 $A+iB$ where $A,B$ a real-valued  $2n\times 2n$ matrices, with $A^t=-A$ and $B^t=B$ and
$\mbox{MSpan$\{j,k\}$}$ will denote the space   of $2n\times 2n$ matrices 
$jC+kD$ where $C,D$ a real-valued  $2n\times 2n$ matrices, with $C^t=C$ and $D^t=D$.\par
  \par
  We are going to show first  an analogous of Theorem \ref{mainth}.
  
  \begin{Th}\label{mainthP} There exists $\varepsilon_0>0$ such that for every $\Gamma \in \dot W^{1,2}(\C,\mbox{MSpan$\{1,i\}$})$ satisfying $\|\Gamma\|_{L^2}\le \varepsilon_0$, every $\Gamma_1 \in \dot W^{1,(2,1)}(\C,{\mathcal{M}}_{2n\times 2n}(\K))$ with  $\|\Gamma_1\|_{L^{2,1}}\le \varepsilon_0$ and every  $G \in L^2(\C,\K)$ solving
 \begin{equation}\label{maineq7}
\partial_L G=(\Gamma+\Gamma_1) G\ ,
\end{equation}
then $G\equiv 0.$ \hfill $\Box$
\end{Th}
As in the case of $2$D codomains the key step to prove Theorem \ref{mainthP} is the following result.
\begin{Prop}\label{L21P}
Let $G \in L^2(\C,\K^{2n})$ be a solution of \eqref{maineq7}. There exists  an $\varepsilon_0>0$  such that
if $\|\Gamma\|_{L^2}\le \varepsilon_0$ and   $\|\Gamma_1\|_{L^{2,1}}\le \varepsilon_0$, then there is a $P\in \dot{W}^{1,2}(\R^2,u(2n,\K))$ and $\chi\in  \dot{W}^{1,(2,1)}(\R^2,\mbox{MSpan$\{1,i\}$})$ such that
$\|\nabla\chi\|_{L^{2,1}}\le  \varepsilon_0$ and
\begin{equation}\label{eqPG}
\partial_{x_1}(PG)-\partial_{x_2}(PiG)=P(-\partial_{x_2}\chi-\partial_{x_1}\chi i+2\Gamma_1)G=2P(-i\partial_L\chi+\Gamma_1) G.
\end{equation}
\end{Prop}
{\bf Proof of Proposition \ref{L21P}.}
If $G$ solves \eqref{maineq7} then as we have seen in \eqref{newsystem2}    for every $P\in u(2n,\K)$ we have
\begin{eqnarray}\label{maineq8}
&&\partial_{x_1}(PG)-\partial_{x_2}(PiG)=P\left[2\Gamma+2\Gamma_1+(P^{-1}\partial_{x_1}P -P^{-1}\partial_{x_2}P i)\right]G
\end{eqnarray}
   
{\bf Step 1.} We introduce the following operator
\begin{eqnarray}\label{estGaugeP}
&& {\mathbf N}\colon \dot W^{1,2}(\C, U(2n,\K))\to  \dot{W}^{-1,2}(\C,\mbox{MSpan$\{1,i\}$})\times L^{2}(\C,\mbox{MSpan$\{j,k\}$})\nonumber\\[5mm]
&&~~~~~P\mapsto(\Pi_{1i}(\partial_{x_1}( {P}^{-1}\partial_{x_1} P)+\partial_{x_2}({P}^{-1}\partial_{x_2}P)),\Pi_{jk}(P^{-1}\partial_{x_1}P-P^{-1}\partial_{x_2} P\,i)\nonumber
\end{eqnarray}
 {\bf Claim 1: }    ${\mathbf N}$ satisfies the following property: 
there is  $\varepsilon_0>0$ and $C>0$ such that for any choice of $V\in  \dot{W}^{-1,2}({\C},\mbox{MSpan$\{1,i\}$})$ and $T\in L^2({\C},\mbox{MSpan$\{j,k\}$})$ satisfying
\begin{equation}\label{smallassbis}
\|V\|_{ \dot{W}^{-1,2}}\le \varepsilon_0,~~\|T\|_{L^{2}}\le \varepsilon_0,
\end{equation}
then
there is $P\in \dot  W^{1,2}(\C,U(2n,\K))$ with
\begin{equation}\label{decompbis}{\mathbf N}(P)=(V,T)\end{equation}
and
\be
\|\nabla P\|_{L^2}\le C(\|V\|_{ \dot{W}^{-1,2}}+\|T|_{L^{2}})\,.
\ee
{\bf Proof of Claim 1.} 
The proof of Claim 1  is very similar to that of Lemma \ref{lm-invertN-a}, therefore we will   sketch  only the main arguments.
For every $P_0\in U(n,\K)$ we consider perturbations of the type: $P=P_0e^{tU}$ where $U\in u(2n,\K)$ and we set
$\tilde{\mathbf N}_{P_0}(U)={\mathbf N}(P_0e^{U})$ and
$$D\tilde{N}_{P_0}(0)= \frac{d}{dt} \tilde{N}_{P_0}(tU) \Big|_{t=0}=: L_{P_0}(U)$$

We have
\begin{equation*}
\label{diffP}
\begin{array}{l}
\ds L_{P_0}(U) :=\left(\Pi_{1i}\left(\Delta {U}+\partial_{x_1}[P_0^{-1}\partial_{x_1}P_0 {U}-{U} P_0^{-1}\partial_{x_1}P_0]+\partial_{x_2}[P_0^{-1}\partial_{x_2}P_0 {U}-{U} P_0^{-1}\partial_{x_2}P_0]\right),\right.\\[5mm]
\ds\quad\quad\quad\left.\Pi_{jk}(\partial_{x_1}{U}-\partial_{x_2}{U} i+[P_0^{-1}\partial_{x_1}P_0 {U}-{U} P_0^{-1}\partial_{x_1}P_0]-[P_0^{-1}\partial_{x_2}P_0 {U}-{U} P_0^{-1}\partial_{x_2}P_0]i)\right)\,.
\end{array}
 \end{equation*}
In the case $P_0=Id$ we get
\begin{eqnarray*}
\label{diffP2}
 \ds L_{Id}(U)& :=&\left(\Pi_{1i}\left(\Delta {U}\right),\Pi_{jk}(\partial_{x_1}{U}-\partial_{x_2}{U} i)\right)\\
&=&\left( \Delta ({U_0+iU_1}), (\partial_{x_1}{U_2}-\partial_{x_2}{U_3}) j+(\partial_{x_1}{U_3}+\partial_{x_2}{U_2}) k\right).
\end{eqnarray*}
Now by arguing exactly  as in the proof of Theorem \ref{mainth} one can prove that
 if $\varepsilon_0$ in \eqref{smallassbis}
  is small    enough  then $L_{P_0}$ with 
$\|\nabla P_0\|_{L^2}<\varepsilon_0$ is invertible,  therefore  {\bf the Claim 1 holds}.\par 
{\bf Step 2.} 
 From Step 1 it follows that  if $\|\Gamma\|_{L^2}<\varepsilon_0$ then there is $P\in \dot  W^{1,2}(\C,U(2n,\K))$ such that
\begin{equation} \label{est10bis}
\left\{\begin{array}{c}
  \Pi_{1i}(\partial_{x_1}(  P ^{-1}\partial_{x_1}   P )+\partial_{x_2}(  P ^{-1}\partial_{x_2}   P ))=0 \\[5mm]
 \Pi_{jk}(  P ^{-1}\partial_{x_1}   P -  P ^{-1}\partial_{x_2}   P \ i)=-2\Gamma\,.
 \end{array}
 \right.
   \end{equation}
 From the first equation in \eqref{est10bis} it follows  the existence of  $\chi
\in \dot W^{1,2}(\C,\mbox{MSpan$\{1,i\}$})$ such that
 \begin{equation}\label{est12bis}
 \left\{\begin{array}{c}
 \Pi_{1i}(P^{-1}\partial_{x_1} P)=-\partial_{x_2}\chi\\[5mm]
 \Pi_{1i}(P^{-1}\partial_{x_2} P)=\partial_{x_1}\chi\,.
 \end{array}\right.
 \end{equation}
 
From \eqref{est12bis} it follows in particular that
\begin{eqnarray}\label{est13-aabis}
-\Delta\chi&=& \partial_{x_2}\left(\Pi_{1i}(  P^{-1}\partial_{x_1}   P )\right)- \partial_{x_1}\left(\Pi_{1i}(  P^{-1}\partial_{x_2}  P )\right)\nonumber\\[5mm]
&=&\Pi_{1i}\left( \partial_{x_2}(  P^{-1}\partial_{x_1}  P )- \partial_{x_1}(  P^{-1}\partial_{x_2}  P )\right)\,.
\end{eqnarray}
The right hand side of \eqref{est13} is a sum of Jacobians, hence it is in the Hardy space ${\mathcal{H}}^1(\C)$.
This implies that $\nabla\chi\in L^{2,1}(\C),$ with $\|\nabla\chi\|_{L^{2,1}}\lesssim \|\nabla  P \|^2_{L^2}.$\par

 We have
 \begin{eqnarray}\label{est13bis}
   P^{-1}\partial_{x_1}  P -  P^{-1}\partial_{x_2}  P \ i&=&\Pi_{1i}(  P^{-1}\partial_{x_1}  P )+\Pi_{jk}(  P^{-1}\partial_{x_1}  P )
 -\Pi_{1i}(  P^{-1}\partial_{x_2}  P )\ i-\Pi_{jk}(  P^{-1}\partial_{x_2}  P )\ i\nonumber\\[5mm]
 &=&-\partial_{x_2}\chi-\partial_{x_1}\chi\ i+\Pi_{jk}(  P^{-1}\partial_{x_1}  P  -P^{-1}\partial_{x_2}  P \ i) \\[5mm]
 &=& -\partial_{x_2}\chi-\partial_{x_1}\chi\ i-2\Gamma\,.\nonumber
 \end{eqnarray}
 
By combining \eqref{maineq8}
\eqref{est12bis}    and \eqref{est13bis} we get
\begin{equation}\label{eqPGbis}
\partial_{x_1}(PG)-\partial_{x_2}(PiG)=P(2\Gamma_1-\partial_{x_2}\chi-\partial_{x_1}\chi i)G=2P(-i\partial_L\chi+\Gamma_1)G
\end{equation}
 and we conclude the proof  of the proposition \ref{L21P}.~\hfill$\Box$
 \par
 \medskip
{\bf Proof of Theorem \ref{mainthP}.}  By arguing as in the end of the proof of Theorem \ref{mainth}   from \eqref{eqPGbis}  we deduce that
\begin{eqnarray}\label{test}
\|PG\|_{L^{2,\infty}}&\lesssim& ( \|\nabla \chi \|_{L^{2,1}}+\|\Gamma_1\|_{L^{2,1}})\|PG\|_{L^{2,\infty}}\nonumber\\[5mm]
&\lesssim&\varepsilon_0\ \|PG\|_{L^{2,\infty}}.
\end{eqnarray}
If $\varepsilon_0$ is small enough then $G\equiv 0.$ This concludes the proof of theorem~\ref{mainthP} and therefore of theorem \ref{bootstrap2}.\hfill $\Box$
\par
\medskip
 From theorem \ref{mainthP} it follows theorem~\ref{th-intro-1} in the general case $n\ge 2.$  The proof is the same of that of Proposition \ref{regul} and therefore we omit it.   }

 \section{Proof of theorem~\ref{th-intro-3-bis}}
A standard covering argument gives that, modulo extraction of a subsequence, there exist finitely many points $a_1\ldots a_Q$ such that, for any $\delta>0$
\begin{eqnarray}\label{eps}
&&\lim_{k\rightarrow +\infty}\inf\lf\{\rho>0\ ;\ \int_{B_\rho(x)}|\nabla S_k|^2(y)\ dy=\frac{\ep_0^2}{2}\ \mbox{ where }\ x\in \R^2\setminus\cup_{i=1}^Q B_\delta(a_i)\rg\}>0\nonumber\\&&
\end{eqnarray}
where $\ep_0>$ is given by the epsilon-regularity theorem ~\ref{bootstrap2}.  Theorem  ~\ref{th-intro-1} implies then that $u_k\rightarrow u_\infty$ strongly
in $L^2_{loc}(\R^2\setminus\{a_1\ldots a_Q\})$ hence we can pass in the limit in the equation away from the points and one gets
\[
\mbox{div}\,(S_\infty\,\nabla u_\infty)=0\quad\quad\mbox{ in }{\mathcal D}'(\R^2\setminus\{a_1\ldots a_Q\})\,.
\]
It remains to establish the point removability. Since $S_\infty\,\nabla u_\infty=\nabla (S_\infty\, u_\infty)-\nabla S_\infty\, u_\infty\in W^{-1,2}+L^1(\R^2)$ a classical result on distributions supported by points gives the existence of $\al_1\ldots\al_Q\in {\R}^n$ such that
\[
\mbox{div}\,(S_\infty\,\nabla u_\infty)=\sum_{i=1}^H\al_i\, \delta_{a_i}\quad\quad\mbox{ in }{\mathcal D}'(\R^2)\,.
\]
We pick a point $a_{i_0}$ arbitrary and we consider an axially symmetric function $\chi$ centered at $a_{i_0}$ such that  $\chi\equiv 1$ in a neighborhood of $a_{i_0}$ and
Supp$\chi\subset B(a_{i_0},r)$ where $0<r<\inf_{i\ne j}|a_i-a_j|$. We have
\[
0=\int_{B(a_{i_0},r)}\nabla \chi\cdot S_k\,\nabla u_k\ dx
\]
Because of the weak convergence of $\nabla u_k$ towards $\nabla u_\infty$ in $L^q$ for any $q<2$ away from the points $a_1\ldots a_Q$ and the strong convergence of $S_k$ towards $S_\infty$ in any $L^p_{loc}$ for $p<+\infty$ we have
\[
0=\int_{B(a_{i_0},r)}\nabla \chi\cdot S_\infty\nabla u_\infty\ dx
\]
which gives $\al_{i_0}=0$. This concludes the proof of theorem~\ref{th-intro-3-bis}.\hfill $\Box$

\end{document}